\input amstex
\input amsppt.sty
\magnification\magstep1

\def\ni\noindent
\def\sbs{\subset}

\def\as{\operatorname{asdim}}
\def\diam{\operatorname{diam}}

\def\ds{\displaystyle}

\def\AN{\operatorname{AN-asdim}}
\def\diam{\operatorname{diam}}
\def\Xbar{\overline{X}}
\def\dim{\operatorname{dim}}

\def\st{\operatorname{st}}
\def\Nerve{\operatorname{Nerve}}
\def\mult{\operatorname{mult}}
\def\cl{\operatorname{cl}}
\def\mesh{\operatorname{mesh}}
\def\Lip{\operatorname{Lip}}
\def\interior{\text{int}}
\def\Kcal{\Cal K}
\def\Fcal{\Cal F}
\def\R{\text{\bf R}}
\def\C{\text{\bf C}}

\def\N{\text{\bf N}}

\def\g{\gamma}
\def\G{\Gamma}
\def\e{\epsilon}
\def\la{\lambda}

\def\sC{\Cal C}

\def\Ecal{\Cal E}

\def\sU{\Cal U}

\def\p{\partial}
\def\x{\times}
\def\Vcal{\Cal V}
\def\sU{\Cal U}
\def\Ucal{\Cal U}

\def\Wcal{\Cal W}

\def\p{\partial}
\def\x{\times}

\def\mk{\vskip .1in}
\def\norm#1{\| #1 \|}

\def\ds{\displaystyle}

\def\Xbar{\overline{X}}
\def\ident#1{\text{id}_{#1}}

\def\Bd{\dot{B}}

\hoffset= 0.0in
\voffset= 0.0in
\hsize=32pc
\vsize=38pc
\baselineskip=24pt
\NoBlackBoxes
\topmatter
\author
A.N. Dranishnikov, J. Smith
\endauthor

\title
On asymptotic Assouad-Nagata dimension
\endtitle
%\date{05.08.99}
%\enddate
\abstract
For a large class of metric spaces $X$ including discrete groups
we prove that the
asymptotic Assouad-Nagata dimension $\AN X$ of $X$
coincides with the covering dimension $\dim(\nu_LX)$ of the
Higson corona of $X$ with respect to the sublinear coarse structure
on $X$. Then we apply this fact to prove the equality
$\AN(X\times\R)=\AN X+1$.
We note that the similar equality for Gromov's asymptotic dimension
$\as$ generally fails to hold [Dr3].

Additionally we construct an injective map $\xi:cone_{\omega}(X)\setminus[x_0]\to\nu_LX$
from the asymptotic cone without the base point to the sublinear Higson corona.
\endabstract

\thanks
The first author was partially supported by an NSF grant. Also
he would like to thank the Max-Planck Institute
f\"ur Mathematik for hospitality.
\endthanks

\address University of Florida, Department of Mathematics, P.O.~Box~118105,
358 Little Hall, Gainesville, FL 32611-8105, USA
\endaddress


\email  dranish\@math.ufl.edu 
justins\@math.ufl.edu
\endemail

\endtopmatter

\document
\head \S1 Introduction \endhead
The Assouad-Nagata dimension was introduced in the 80s by Assouad [As1],[As2]
under the name Nagata dimension. Recently this notion was revived in
the asymptotic geometry due to works of Lang and Schlichenmaier [LSch], and
Buyalo and Lebedeva [Bu], [BL].
The concept takes into account the dimension of
a metric space on all scales. In this paper we consider only the large scale
version of it. Note that the asymptotic version of the Assouad-Nagata dimension
agrees with the original for our main source of examples of metric spaces -
finitely generated discrete groups with the word metric.
Like in the case of Gromov's asymptotic dimension, the Assouad-Nagata dimension
is a group invariant.

A certain analogy between the asymptotic Assouad-Nagata
dimension $\AN$ and the asymptotic dimension $\as$ invites one to transfer
the asymptotic dimension theory [Gr],[Dr1],[Dr2],[Dr3],[DKU],[BD1],[BD2],[BD3],
[DZ],[Ro2] to the asymptotic Assouad-Nagata dimension. It was partially done
in [LSch], [BDHM], [BDLM]. Namely, the theorem on
embedding into a product of trees,
the characterization of $\as$ in terms of map extension, the union
theorems, and the Hurewicz type theorem  were successfully extended
to the case of the Assouad-Nagata dimension. In this paper we extend to
the Assouad-Nagata dimension the theorem that characterizes the asymptotic
dimension of a metric space $\as X$ as the
covering dimension of the Higson corona $\dim\nu X$ ([Dr1],[DKU]).
For that we introduce a coarse structure $\Ecal_L$ on a metric space $X$ called
the sublinear coarse structure and show that the covering dimension of
the Higson corona $\nu_LX$ of this coarse structure is exactly the
asymptotic Assouad-Nagata dimension of the space, provided the latter is finite.
Contrary to the case of the classic Higson corona, the sublinear
Higson corona $\nu_L$ behaves nicely under the product
with reals. Namely, there is a decomposition:
$\nu_L(X\times\R)=\nu_LX\times(-1,1)\cup\nu_L\R$. In particular,
it implies $\dim\nu_L(X\times\R)=\dim\nu_LX+1$. We prove that $\AN X=\dim\nu_LX$
for sufficiently symmetric spaces $X$ like discrete groups (it is not true in
general). Then for such spaces $\AN(X\times\R)=\AN X+1$. This is an analog of the
classical Morita formula from the dimension theory: $\dim(X\times\R)=\dim X+1$.
We note that the Morita formula generally does not hold for asymptotic
dimension [Dr3].

\mk

{\bf Coarse structures.} A coarse structure $\sC$ on a set $X$ is a family
of subsets $E\subset X\times X$ that contains the
diagonal $\Delta_X$ and is closed taking
finite unions, subsets, inverses, and compositions.  The elements
of $\sC$ are called {\it controlled} sets (see [HR],[Ro2],[DH]).

Suppose that $X$ is a topological space. Then
a set $E\subset X\times X$ is called {\it proper} if both $E[K]$ and
$E^{-1}[K]$ are relatively compact for a relatively compact set $K\subset X$,
where $E[K]$ is the set of all $x'$ such that there is $x\in K$ with $(x',x)\in
E$. We use the notations $E_x=E[x]$ and $E^x=E^{-1}[x]$ for $x\in X$.

A subset $B\subset X$ of a coarse space is {\it bounded} if $B\times B$ is
controlled. A map between coarse spaces $f:(X,\sC)\to(x',\sC')$ is called
a {\it proper} if the preimage $f^{-1}(B)$ of every bounded set is bounded.  
A map between coarse spaces $f:(X,\sC)\to(x',\sC')$ is called a {\it  coarse
morphism} if it is coarsely proper and $(f\times f)$ takes controlled sets to
controlled.

Suppose that $X$ is a topological space.  We say
that a coarse space $(X, \Ecal)$ is
{\it consistent} with the topology on $X$ if $B \subset X$ is (coarsely) bounded
if and only if $B$ is relatively compact (i.e., bounded sets coincide
with relatively compact sets).  One can easily show a consistent coarse space
$(X, \Ecal)$ is coarsely connected and each $E \in \Ecal$ is proper.  If
$X$ is a locally compact Hausdorff topological space, we
say that $(X, \Ecal)$ is {\it proper} if $(X, \Ecal)$ is consistent with the
topology and if $\Ecal$ contains a neighborhood of the diagonal.

\mk

{\bf Compactifications.} Let $\bar{X}$ be a compactification of a locally
compact space $X$, and let
$V$ be an open subset of $X$.  Then there is a unique maximal open set
$\widetilde{V}$ in $\bar{X}$ such that $\widetilde{V} \cap X = V$.
In fact,
$\widetilde{V} = \bar{X} \setminus \overline{X \setminus V}$.  One can
show that
$\widetilde{V} \subset \overline{V}$.

The following propositions are obvious.

\proclaim{Proposition 1.1} Let $\bar{X}$ be a compactification of a locally
compact space
$X$, and let $\nu X = \bar{X} \setminus X$.  Then
$\{ \widetilde{V} \cap \nu X : V \text{ \ is open in \ } X \}$ forms a
basis for $\nu X$.
\endproclaim

\proclaim{Proposition 1.2} Let $\bar{X}$ be a compactification of a locally compact space
$X$, and let $\nu X = \bar{X} \setminus X$.  Suppose $U \subset X$ is open and suppose
$x \in \widetilde{U} \cap \nu X$.  Then there is a set $V \subset U$ open in $X$ such
that $x \in \widetilde{V} \cap \nu X$ and $\bar{V} \subset \widetilde{U}$.
\endproclaim

\demo{Proof} Let $W$ be an open subset of $\Xbar$ such that
$x \in W \subset \overline{W} \subset \widetilde{U}$, and
set $V = W \cap X$.  We have that $V$ is open in $X$, $V \subset W$, and
$W \subset \widetilde{V}$ by definition.  Thus,
$x \in \widetilde{V} \cap \nu X$ and
$\overline{V} \subset \overline{W} \subset \widetilde{U}$.  This completes
the proof.
\qed
\enddemo

\proclaim{Proposition 1.3} Let $\bar{X}$ be a compactification of a locally compact space
$X$, and let $\nu X = \bar{X} \setminus X$.  Suppose $U$ is an open subset of
$\nu X$ and $x \in U$.  Then there is a set $V$ which is open in $X$,
$x \in \widetilde{V} \cap \nu X$, and $\bar{V} \cap \nu X \subset U$.
\endproclaim

\demo{Proof} Choose $W_1$ open in $\overline{X}$ such that
$U = W_1 \cap \nu X$ and take
$W_2$ open in $\overline{X}$ such that
$x \in W_2 \subset \overline{W_2} \subset W_1.$  Set
$V = W_2 \cap X$.  Thus, $W_2 \subset \widetilde{V}$, hence
$x \in W_2 \cap \nu X \subset \widetilde{V} \cap \nu X$.  Also,
$\overline{V} \cap \nu X \subset \overline{W_2} \cap \nu X \subset W_1 \cap \nu X = U$
since $V \subset W_2$. \qed
\enddemo

Suppose $\overline{X}$ is a compactification of the locally compact Hausdorff space $X$.  Then
$(X, \overline{X})$ will be called a compactified pair.  Now suppose,in addition, that
$\Ecal$ is a coarse structure which is
consistent with the topology on $X$.  We say that $f: X \to \C$ is a {\it Higson function},
denoted $f \in C_h (X, \Ecal)$, if for every $E \in \Ecal$ and
every $\e > 0$, there is a compact set
$K$ such that
$|f(x) - f(y)| < \e$ whenever $(y,x) \in E \setminus K \times K$.
Then by the GNS theorem there is a compactification $h_{\Ecal}X$ of
$X$ called the {\it Higson compactification} such that the algebra of Higson functions
$C_h (X, \Ecal)$ is isomorphic to
$C( h_{\Ecal} X)$.  We define $h(X, \Ecal) = (X, h_{\Ecal}X)$.  The
{\it Higson corona} is defined by $\nu_{\Ecal} X = h_{\Ecal} X \setminus X$.

On the other hand, suppose we have a compactified pair $(X, \overline{X})$.  Let
$\Ecal_{\overline{X}}$ be those $E \subset X \times X$ for which
$\overline E \setminus X\times X \subset\Delta_{\p X}$,
where
$\partial X = \overline{X} \setminus X$,  $\overline{E}$ denotes
the closure of $E$ in $\overline{X} \times \overline{X}$, and $\Delta_A$ denotes
the diagonal in $A\times A$.  Then
$\Ecal_{\overline{X}}$ is a coarse structure on $X$
which is consistent with the topology on $X$.  We will sometimes use the notation of
Roe [Ro2],
$t \overline{X}$, instead of $(X, \Ecal_{\overline{X}})$.

The following generalizes a definition from \cite{DKU}.

DEFINITION.  For a general coarse space $(X, \Ecal)$,
a finite system $E_1, \ldots, E_n$ of subsets of $X$ {\it diverges} if
$$ \bigcap_{i = 1}^{n} F[E_i]$$
is bounded for each $F \in \Ecal$.

\proclaim{Theorem 1.4} Let $X$ be a locally compact Hausdorff space equipped
with a coarse structure $\Ecal$ that is consistent with the topology.
For a finite system $E_1, \ldots, E_n$ of
subsets of $X$, if $\nu X \cap [ \cap_{i =1}^{n} \overline{E_i} ] = \varnothing$,
then the system $E_1, \ldots, E_n$ diverges.
\endproclaim

\demo{Proof} We let $\overline{X}$ denote the Higson compactification with respect
to this coarse structure.  Suppose that the system $E_1, E_2, \ldots, E_n$
does not diverge; so there is a controlled set $F$ such that
$\cap_{i=1}^n F[E_i]$ is not bounded.  Thus, for each compact
subset $K$ of $X$, there is an $x_K \in (\cap_{i=1}^n F[E_i]) \setminus K$.
The collection of compact subsets of $X$, ordered by inclusion,
forms a directed set.  We denote it by $\Kcal$.
In particular, $\{ x_K \}_K$ is a net.  Since $\overline{X}$ is
compact, there is a convergent subnet $\{ x_{g(\lambda)} \}_{\lambda}$
(here, $g: \Lambda \rightarrow \Kcal$ is an order-preserving map between
directed sets such that $g(\Lambda)$ is cofinal in $\Kcal$).  Thus, for
$K$ compact, we have by cofinality that there is a $\la_0$ such that
$g(\lambda) \supset K$ whenever $\la \geq \la_0$.  Hence, for
$\la \geq \la_0$, we
have $x_{g(\la)} \in X \setminus g(\la) \subset X \setminus K$.  As
$X$ is locally compact, this means
that $x := \lim_{\la} x_{g(\la)} \in \nu X.$

Now fix $1\leq i \leq n$.  Then $x_{g(\la)} \in F[E_i]$ for
all $\la$; for each $\la$, choose $y_{\la} \in E_i$ such that
$(x_{g(\la)}, y_{\la}) \in F$.  Since $F \in \Ecal$, by Proposition 2.45 (a)
of \cite{Ro2}, we have $F \in t\overline{X}$.  Thus,
since $x_{g(\la)} \rightarrow x \in \nu X$, we must have
$y_{\la} \rightarrow x$.  Thus, $x \in \overline{E}_i$
for each $i$ and
so $\nu X \cap [ \cap_{i =1}^{n} \overline{E_i} ] \neq \varnothing$.
\qed
\enddemo

The following theorem can be found in [Ro2].

\proclaim{Theorem 1.5} Let $X$ and $Y$ be locally compact, Hausdorff spaces
equipped with coarse structures $\Ecal$ and $\Fcal$, respectively, which
are consistent with the topologies.  If $f:X \rightarrow Y$ is a
coarse, continuous map, then $f$ extends to
a continuous map $\overline{f}:h_{\Ecal}X \rightarrow h_{\Fcal}Y$
such that $\overline{f}(\nu X) \subset \nu Y$.  \endproclaim

In fact, writing $h(f)$ rather than $\overline{f}$, we have
that $h$ is a functor from the category of consistent coarse structures
(with the topological space being locally compact) with continuous coarse maps
as morphisms, to the category of compactified pairs with continuous maps
preserving boundary as morphisms.  Also, there is a functor
$\nu$ from the category of proper coarse structures with coarse maps as morphisms,
to the category of compact spaces with continuous maps as morphisms; it sends
$(X, \Ecal)$ to $\nu_{\Ecal} X$.  For the latter, see proposition 2.41
of \cite{Ro2}.

{\bf Asymptotic dimension.} We recall Gromov's definition of
{\it asymptotic dimension} of a metric space [Gr].

DEFINITION. The asymptotic dimension of a metric space $X$
does not exceed $n$, $\as X\le n$, if for every $r>0$ there are
$r$-disjoint uniformly bounded families $\sU^0,\dots,\sU^n$ of subsets of $X$
such that the union $\cup\sU^i$ is a cover of $X$.

There are equivalent reformulations [Gr],[BD2]:
\proclaim{Proposition 1.6} Let $(X,d)$ be a metric space.  The following are
equivalent:
\roster
\item $\as{X} \leq n$;
\item  for
every $\e > 0$, there is $b>0$ and an
$\e-$Lipschitz, $b$-cobounded map $p:X \rightarrow P$ to an
$n-$dimensional uniform simplicial complex $P$;
\item for all $r \geq 0$, there is a uniformly bounded cover $\Ucal$ of $X$ such that
the Lebesgue number $L(\Ucal )\geq r$ and $\Ucal$ has multiplicity $\leq n+1$.
\endroster
\endproclaim
Here a map of a metric space to a simplicial complex $f:X\to P$ is called
$b$-cobounded if $\diam(f^{-1}(\sigma))\le b$ for every simplex $\sigma\subset
P$. A uniform metric on a simplicial complex $P$ is the restriction of the
Euclidean metric from $\ell_2(P^{(0)})$, the Hilbert space spanned by
the vertices of $P$, to $P\subset\ell_2(P^{(0)})$.  By the multiplicity 
of a cover $\Ucal$ of $X$ we mean the minimum number $m$ such that every 
intersection of $m+1$ distinct elements of $\Ucal$ is empty.  We will sometimes 
denote this number by $\mult\Ucal$.

DEFINITION.  The {\it asymptotic Assouad-Nagata dimension} of a
metric space $X$
does not exceed $n$, $\AN{X} \leq n$, if there is a $c>0$ and an
$r_0 > 0$ such that
for every $r \geq r_0$, there is a cover $\Ucal$ of $X$ such that
$\mesh{\Ucal} \leq cr$,
$L(\Ucal) > r$, and $\Ucal$ has multiplicity $\leq n+1$.

This has many aliases, including {\it asymptotic dimension of linear type} and
{\it asymptotic dimension with Higson property}.  For discrete metric spaces,
in particular for discrete finitely generated groups, this definition coincides
with the Assouad-Nagata dimension.

There are analogous reformulations:

\proclaim{Proposition 1.7} Let $(X,d)$ be a metric space.  The following are
equivalent:
\roster
\item $\AN{X} \leq n$;
\item there is a $C>0$ and an $\e_0 >0$ such that for
all $ \e \leq \e_0$ ($\e > 0$), there is an
$\e-$Lipschitz, $C/ \e -$cobounded map $p:X \rightarrow P$ to an
$n-$dimensional simplicial complex $P$;
\item there is a $C>0$ and an $r_0 > 0$ such that, for all
$r \geq r_0$, there are $r-$disjoint families $\Ucal_0, \Ucal_1, \ldots, \Ucal_n$
(of subsets of $X$) such that $\Ucal = \cup_i \Ucal_i$ is a cover of $X$ and
$\mesh \Ucal \leq Cr$;
\item there is a $C>0$ and an $r_0 >0$ such that for all $r \geq r_0$, there is 
a cover $\Ucal$ of $X$ such that
$\mesh \Ucal \leq Cr$ and $B_r(x)$ meets at most $n+1$ elements
of $\Ucal$ for each $x \in X$.
\endroster
\endproclaim

In (1),(3), and (4), we can take the covers to be open.  

\head \S2 The sublinear coarse structure\endhead

We consider a proper metric
space $(X, d)$ with basepoint $x_0$ and define
$\norm{x} = d(x, x_0)$.  We will sometimes write $B_r$ to indicate
$B_r(x_0)$, the open ball of radius $r$ centered at $x_0$.

DEFINITION. We define the {\it sublinear coarse structure}, denoted $\Ecal_{L}$,
on $X$ as follows:
$$ \Ecal_{L} = \{ E \subset X\times X : E \text{ proper},
 \lim_{x \rightarrow \infty} \frac{\sup_{y \in E_x} d(y,x) }{\norm{x}} = 0 =
\lim_{x \rightarrow \infty} \frac{\sup_{y \in E^x} d(x,y) }{\norm{x}} \}. $$
By the statement
$\lim_{x \rightarrow \infty} \frac{\sup_{y \in E_x} d(y,x) }{\norm{x}} = 0$, we
mean that for each $\e > 0$, there is a compact subset $K$ of $X$ containing $x_0$
(equivalently, an $r\geq 0$) such that
$$ \frac{\sup_{y \in E_x} d(y,x) }{\norm{x}} \leq \e $$
for all $x \notin K$ (respectively, for all $x$ with $\norm{x} > r$).  It would perhaps
be better to think of this as $\lim_{\norm{x} \rightarrow \infty}$.  In the 
event that $E_x = \varnothing$, we define $\sup_{y \in E_x}d(y,x) = 0$.  
We leave to the reader to check that $\Ecal_{L}$ is indeed a coarse structure
and that it does not depend on the choice of basepoint.
The Higson corona for the sublinear coarse structure on $X$ will be denoted
by $\nu_LX$.  We will sometimes call this corona the sublinear Higson corona, 
to eliminate possible confusion with the usual Higson corona.  

We recall that a mapping $f:X \to Y$ between metric spaces
is said to be a {\it quasi-isometry} if there are numbers $\la> 0$, $C\geq 0$,
and $D \geq 0$ such that
$$  \frac{1}{\la} d(x,y) - C \leq d(f(x),f(y)) \leq \la d(x,y) + C  $$
and every point of $Y$ is within distance $D$ of $\phi(X)$.

\proclaim{Proposition 2.1} Let $X$ and $Y$ be proper metric spaces.  If $f:X \rightarrow Y$ is
a quasi-isometry, then it is a coarse equivalence with respect to the sublinear coarse structures.
\endproclaim

\demo{Proof}  Fix a basepoint $x_0 \in X$; set $y_0 = f(x_0) \in Y$.
Choose $\la > 0$ and $C \geq 0$ such that
$\frac{1}{\la} d(x,y) - C \leq d(f(x),f(y)) \leq \la d(x,y) + C$.  Let
$g$ be a quasi-isometry such that $f\circ g$ and $g\circ f$ are close to
the respective identity functions.
It is clear that $f$ is proper and the image under $f$ of a bounded set is bounded.
Let $E$ be a controlled set in the sublinear coarse structure on $X$. 
It is not hard to show that $f\x f (E)$ is proper.  

Let $\e > 0$ be given.  There is a bounded $K' \subset X$ such that
$\frac{ \sup_{y \in E_x} d(y,x) }{\norm{x}} \leq \frac{\e}{4\la^2}$ whenever
$x \notin K'$.  Set $K = B(x_0, 2\la C) \cup B(x_0, \frac{4 \la C}{\e}) \cup K'$.  Note
that $f(K)$ is bounded.
Suppose that $z \notin f(K)$.  If $(f \times f)(E)_z = \varnothing$, then 
$\ds \frac{ \sup{ \{ d(z',z): z' \in (f \times f)(E)_z \} }}{\norm{z}} = 0
< \e$ and we are finished.
Now assume that $z' \in (f \times f)(E)_z$; so $z' = f(x')$ and $z = f(x)$
for some $x,x' \in X$ with $(x', x) \in E$.  We have that $x \notin K$.

Then
$$ \norm{f(x)} \geq \frac{1}{\la} \norm{x} - C = \frac{2\norm{x} - 2\la C}{2 \la }
= \frac{ \norm{x}}{2 \la} + \frac{ \norm{x} - 2 \la C}{ 2 \la}
\geq \frac{\norm{x}} {2 \la}$$
and so
$$ \frac{d(z',z)}{\norm{z}} = \frac{d(f(x'),f(x))}{\norm{f(x)}}
\leq 2 \la \frac{\la d(x',x) + C}{\norm{x}}
= 2\la^2 \frac{d(x',x)}{\norm{x}} + \frac{2\la C}{\norm{x}}
\leq \e. $$
Thus, we have
$\ds \frac{ \sup{ \{ d(z',z) : z' \in (f \times f)(E)_z \} }}{\norm{z}} \leq \e$.
Since
$\e > 0$ was arbitrary, $(f \times f)(E)$ is controlled and so $f$ is a coarse map.

Similarly, since $g$ is a quasi-isometry, it is a coarse map as well.
Finally, it is
clear that $f \circ g$ and $g \circ f$ are close to the corresponding
identities when $X$ and
$Y$ are equipped with the sublinear coarse structures.
 Thus, $f$ is a coarse equivalence.\qed  \enddemo
\proclaim{Corollary 2.2}
The sublinear coarse structure $\Ecal_{L}$ is well-defined
on finitely generated groups, i.e., the sublinear coarse structure 
for a given group $\Gamma$ is independent of the 
choices of the finite generating set and the basepoint.
In particular, the asymptotic dimension $\as(\Gamma,\Ecal_L)$
associated with $\Ecal_{L}$ is
a group invariant for finitely generated groups.
\endproclaim

Next we give a characterization of divergent systems for the sublinear coarse structure.

\proclaim{Lemma 2.3} Let $(X,d)$ be a proper metric space with basepoint $x_0$,
endowed with the sublinear
coarse structure $\Ecal_{L}$.  For a finite system $E_1, \ldots, E_n$ of
subsets of $X$, the following are equivalent.
\roster
\item $\nu_L X \cap [ \cap_{i =1}^{n} \overline{E_i} ] = \varnothing$;
\item the system $E_1, \ldots, E_n$ diverges;
\item there exist $c, r_0 > 0$ such that
$\max_{1 \leq i \leq n} d(x, E_i) \geq c\norm{x}$ whenever $\norm{x} \geq r_0$.
\endroster
\endproclaim

\demo{Proof} That (1) implies (2) was shown earlier.  We now prove that (2) implies (3). Assuming
that (3) does not hold, then if we let $m$ be a positive integer, and
if we set $c = \frac{1}{4m}$ and $r_0= 2m$, then
there is an $x_m$ such that $\norm{x_m} \geq 2m$ yet
$$ \max_{1 \leq i \leq n} d(x_m, E_i) < \frac{1}{4m} \norm{x_m}. $$
Thus, for each $i$, we have that $d(x_m, E_i) < \frac{1}{4m} \norm{x_m}$, and so
there is an $a^i_m \in E_i$ such that $d(x_m, a^i_m) < \frac{1}{4m}\norm{x_m}$.
We have $d(a^i_m, a^j_m ) < \frac{1}{2m}\norm{x_m}$.  Also,
$$\norm{x_m} \leq \norm{a^i_m} + d(a^i_m, x_m) < \norm{a^i_m} + %
\frac{1}{4m}\norm{x_m}, $$
and hence $(1-\frac{1}{4m})\norm{x_m} < \norm{a^i_m}$ (all $i$).
Since $\frac{1}{4m}<1/2$,
we have $\norm{a^i_m} > \frac{1}{2} \norm{x_m}$.  Thus,
$$ d(a^i_m, a^j_m) < \frac{1}{2m}\norm{x_m} < \frac{1}{m} \norm{a^j_m}  \quad \text{and} \quad %
\norm{a^i_m} > \frac{1}{2} \norm{x_m} \geq m  $$
for all $1 \leq i,j \leq n$.
Take $F_{i,j} = \{ (a^i_m,a^j_m) : m=1,2,\ldots \}$ for each $1\leq i,j \leq n$.
Fixing $i,j$, we temporarily set $G = F_{i,j}$ for convenience, and show that $G$
is controlled.  Since $\norm{a^i_m} \rightarrow \infty$
and $\norm{a^j_m} \rightarrow \infty$ as $m \rightarrow \infty$, it follows that
$G$ is proper.  Now let $\e > 0$ be given, and take $M$ to be a positive
integer for which $1/M < \e$; set $K = \{ a^j_m : 1\leq m < M \}$.
Suppose that $x \notin K$.  If $G_x = \varnothing$, then by our convention we have
$\frac{\sup_{y \in G_x} d(y,x)}{\norm{x}} = 0$.  If $y \in G_x$, then
there is a positive integer $m$ such that $(y,x) = (a^i_{m},a^j_{m})$, and
since $a^j_m = x \notin K$, we must have $m\geq M$; it
follows that
$\frac{d(y,x)}{\norm{x}} = \frac{d(a^i_m, a^j_m)}{\norm{a^j_m}} < \frac{1}{m} < \e$
and so $\frac{\sup_{y \in G_x} d(y,x)}{\norm{x}} \leq \e$.  Thus,
$\lim_{x \rightarrow \infty}\frac{\sup_{y \in G_x} d(y,x)}{\norm{x}} = 0$.
Similarly, $\lim_{x \rightarrow \infty}\frac{\sup_{y \in G^x} d(x,y)}{\norm{x}} = 0$, and
hence $G=F_{i,j}$ is controlled.

Define $F = \cup_{1\leq j \leq n} F_{1,j}$ and $A = \{ a^1_m : m = 1,2,\ldots \}$.  Note
that $A$ is not bounded and $F$ is controlled.  Also,
$F[E_j] \supset F_{1,j}[E_j] \supset A$ for all $1\leq j \leq n$ and hence
$$ \cap_{j=1}^n F[E_j] \supset A, $$
which means that $ \cap_{j=1}^n F[E_j]$ is not bounded.  So (2) does not hold.

It remains to show that (3) implies (1).
Define $F_i = E_i \setminus B_{r_0 + c r_0}$ for $1 \leq i \leq n$.
Let $f: X \rightarrow \R$ be
defined by $f(x) = \sum_{i=1}^{n} d(x, F_i)$.  Note that
$f(x) \geq c \norm{x}$ when $\norm{x} \geq r_0$ since $d(x,F_i) \geq d(x,E_i)$.  Also,
$f(x) \geq c r_0$ when $\norm{x} \leq r_0$; in particular,
$f(x)\geq c \norm{x}$ for all $x$ and $f(x) > 0$ for all $x$.  Define
$g_i : X \to \R$ by $g_i(x) = d(x, F_i ) / f(x)$.

Let $E$ be a controlled set.  Since
$$ \align
|g_i (y) - g_i (x)| & \leq  d(y, F_i)|\frac{1}{f(y)} - \frac{1}{f(x)}| + %
| \frac{d(y, F_i) - d(x, F_i)}{f(x)}| \leq  \frac{d(y, F_i)}{f(x)f(y)} |f(x) - f(y)| \\ 
 & + \frac{d(y,x)}{f(x)} \leq \frac{n d(y,x)}{f(x)} + \frac{d(x, y)}{f(x)} %
 \leq (n+1) \frac{d(x,y)}{c \norm{x}}, 
\endalign $$
we have $\sup_{y\in E_x} |g_i(x)-g_i(y)| \to 0$ as
$x \rightarrow \infty$.  Since $E$ was an arbitrary controlled set, 
$g_i$ (viewed as a map to $\C$) is a Higson function for each
$i$.  Let $G_i : \overline{X} \rightarrow \C$
be the extension of $g_i$ to the Higson compactification.  Since
$\sum_i g_i = 1$, it is immediate that
$\sum_i G_i =1 $ throughout $\Xbar$.  Also, $\overline{F_i} \subset G_i^{-1}(0)$ 
and it is not hard to see that
$\nu X \cap \overline{F_i} = \nu X \cap \overline{E_i}$.
Thus,
$$\nu X \cap ( \cap_{i=1}^n \overline{E_i} ) = %
\nu X \cap ( \cap_{i=1}^n \overline{F_i} ) \subset %
\nu X \cap ( \cap_{i=1}^n G_i^{-1}(0) ) = \varnothing $$
since $\sum_i G_i = 1$ on $\nu X$. \qed
\enddemo

In the case that $n = 2$, we can add another condition.

\proclaim{Lemma 2.4} Let $A$ and $B$ be subsets of a metric space $X$.  Let $x_0 \in X$,
and define $\norm{\cdot}$ as usual.  Also, take $B_r = B_r (x_0)$.  Then
the following are equivalent.
\roster
\item There exist $C, r_0 > 0$ such that $\max{ \{ d(x, A) , d(x, B) \} } \geq C \norm{x}$ whenever
$\norm{x} \geq r_0$;
\item there exist $D, r_1 > 0$ such that $d(A \setminus B_r, B \setminus B_r) \geq Dr$ whenever
$r \geq r_1$.
\endroster
\endproclaim

\demo{Proof} We show (1) implies (2).  Given $C$ and $r_0$, take
$D = C$ and $r_1 = r_0$.  Let $r \geq r_1$, $a \in A \setminus B_r$, and
$b \in B \setminus B_r$.  So $\norm{a} \geq r \geq r_0$.  Thus,
$$Cr \leq C\norm{a} \leq \max{ \{ d(a,A), d(a,B) \} } = \max{\{ 0, d(a,B) \} } = d(a,B) \leq d(a,b)$$
by (1).  So $d(A \setminus B_r, B \setminus B_r) \geq Dr$.

We show (2) implies (1).  Let $D, r_1$ be positive numbers satisfying (2).
Set $r_0 = 2 r_1$ and take $C$ to be a positive number satisfying
$C < \min{ \{ 1/2, D/4 \} }$.  Now let $x \in X$ be such that $\norm{x} \geq r_0 = 2r_1$.
To get a contradiction, suppose that $\max{ \{ d(x,A), d(x,B) \} } < C\norm{x}$.
Then $d(x,A) < C \norm{x}$ and $d(x,B) < C \norm{x}$.  Thus,
there exist $a \in A$ and $b \in B$ such that
$d(x,a) < C\norm{x}$ and $d(x,b) < C\norm{x}$.  So $d(a,b) < 2C\norm{x}$.  We then have
$$ \norm{a} \geq \norm{x} - d(x,a) > \norm{x} - C \norm{x} = (1-C)\norm{x} \geq \norm{x}/2.$$
Similarly, $\norm{b} \geq \norm{x}/2$.  Since $\norm{x}/2 \geq r_1$, we have
by (2) that
$$ \frac{D \norm{x}}{2} \leq d(A \setminus B_{\norm{x}/2}, B \setminus B_{\norm{x}/2}) %
\leq d(a,b)< 2C\norm{x} \leq \frac{D\norm{x}}{2}, $$
a contradiction.  Therefore, $\max{ \{ d(x,A), d(x,B) \} } \geq C\norm{x}$ 
when $\norm{x}\geq r_0$. \qed \enddemo

DEFINITION.  Let $(X, d)$ be a metric space, and let $\Vcal$ be a 
family of open subsets of $X$.  We define the {\it Lebesgue function } associated
with the cover $\Vcal$, denoted $L^{\Vcal}$, by
$$L^{\Vcal} (x) = \sup_{V \in \Vcal} d(x, X \setminus V ).  $$

DEFINITION.  For a proper metric space $(X, d)$ with basepoint $x_0$, we say a function
$f: X \rightarrow [0, \infty )$ is (eventually) {\it at least linear} if
there exist $c, r_0 > 0$ such that $f(x) \geq c \norm{x}$ whenever
$\norm{x} \geq r_0$.

\proclaim{Corollary 2.5} Let $(X,d)$ be a proper metric space endowed with the
coarse structure $\Ecal_L$.  Let $\alpha = \{ O_1, \ldots, O_n \}$ be a
finite family of open subsets of $X$.  Then
$\widetilde{\alpha} = \{ \widetilde{O}_1, \ldots, \widetilde{O}_n \}$ covers the corona
$\nu_L X$ if and only if the Lebesgue function $L^{\alpha}$ is at least linear.
\endproclaim

\demo{Proof} $\widetilde{\alpha} = \{ \widetilde{O}_1, \ldots, \widetilde{O}_n \}$ covers the corona
$\nu X$ iff
$\nu X \setminus (\cup_i \widetilde{O}_i ) = \varnothing$, iff
$\nu X \cap ( \Xbar \setminus \cup_i \widetilde{O}_i ) = \varnothing $, iff
$\nu X \cap (\cap_i (\Xbar \setminus \widetilde{O}_i)) = \varnothing$, iff
$\nu X \cap (\cap_i \overline{X \setminus O_i} ) = \varnothing$ by the comments preceding
Proposition 1.1, iff
the system $X \setminus O_1, \ldots, X \setminus O_n$ diverges, which, by lemma 2.3
above, is
true if and only if $L^{\alpha}$ is at least linear.
\qed
\enddemo

\proclaim{Corollary 2.6} Let $(X, d)$ be a proper metric space, and
let $A$ be a closed subspace of $X$ equipped with the restricted metric.  Then
the embedding $A \rightarrow X$ extends to an embedding $h_L A \rightarrow h_L X$ on the
compactifications and
induces an embedding $\nu_L A \rightarrow \nu_L X$ on the coronas.
\endproclaim

\demo{Proof} Let $x_0 \in A$ be the basepoint for both $A$ and $X$, and write
$hX = h_L X$ and $\nu X = \nu_L X$.
Since the inclusion map $i: A \rightarrow X$ is continuous and coarse,
$i$ extends to a continuous map from $h A$ to $h X$ such that
$i(\nu A) \subset \nu X$.  To
prove the result, it suffices to show that $i$ is injective on $\nu A$.  Let
$x_1,x_2 \in \nu A$ with $x_1 \neq x_2$.  First, we can find
disjoint open subsets $U_1,U_2$ of $\nu A$
such that $x_j \in U_j$ for $j=1,2$.  Applying Proposition 1.3,
there are open subsets $V_1$ and $V_2$ of $A$ such that $x_j \in \widetilde{V}_j \cap \nu A$
and $\cl_{hA} {V_j}\cap \nu X \subset U_j$ for $j=1,2$.  So
$\nu A \cap \cl_{hA} {V_1} \cap \cl_{hA} {V_2} \subset U_1 \cap U_2 = \varnothing$,
and by Lemmas 2.3 and 2.4, we have
that there exist $c,r_0 > 0$ such that
$$ d(V_1 \setminus B_r , V_2 \setminus B_r ) =
d|_{A} (V_1 \setminus B_r , V_2 \setminus B_r  ) \geq cr$$
whenever $r \geq r_0$.  Thus,
$\nu X \cap \cl_{hX} {V_1} \cap \cl_{hX} {V_2} = \varnothing$.  But
$i(x_j) \in i( \cl_{hA}{V_j}) \cap \nu X = \cl_{hX}V_j \cap \nu X$ for $j=1,2$,
which means that $i(x_1) \neq i(x_2)$. \qed   \enddemo

\mk

{\bf Algebra of functions.}
We define a subalgebra $U(X)=U(X,x_0)$ of $C(X)$ as follows: $f: X \rightarrow \C$ is in
$U(X)$ if and only if $f$ is bounded, continuous, and there exists a
$c = c_f$ such that
$$ | f(x)-f(y) | \norm{x} \leq c d(x,y). $$

It is not difficult to check that $U(X)$ is closed under addition, multiplication,
and complex conjugation.

REMARK. In the definition, the continuity condition is almost
unnecessary.  It is not hard to see that the property
$|f(x)-f(y)| \norm{x} \leq c d(x,y)$ implies that $f$ is continuous for
$x \neq x_0$.

It is easy to show that $U(X)$ separates points and closed sets.
It is also clear that $U(X)$ is, in general, not complete.  We set
$C'(X) = \overline{U(X)}$, where the bar represents
closure in $C(X)$ with the uniform metric.  
So $C'(X)$ is a $C^{*}-$algebra which separates points and closed sets.  Thus,
by the GNS Theorem
we can extract a compactification of $X$ which will be called the {\it sublinear
compactification}:
\proclaim{Proposition 2.7} With the notation above, there is a compactification
$\Xbar$ of $X$ such that $C'(X) = C(\Xbar)$.  \endproclaim

Let $h_LX$ be the Higson compactification for $\Ecal_L$, the sublinear
coarse structure on $X$.  We have
the following.

\proclaim{Proposition 2.8} We have $C'(X) \subset C_h(X, \Ecal_L)$, and hence there
is a surjective, continuous map $h_L X \rightarrow \Xbar$ which
extends the identity.  \endproclaim
\demo{Proof} Let $f \in U(X,x_0)$, and let $c$ be a constant such that
$|f(x)-f(y)| \norm{x} \leq c d(x,y)$.  Let $E \in \Ecal_L$.  So
$$\lim_{x \rightarrow \infty} \sup_{y \in E_x} |f(x)-f(y)| \leq %
c \lim_{x \rightarrow \infty} \frac{\sup_{y \in E_x} d(y,x)}{\norm{x}} = 0.  $$
But $E \in \Ecal_L$ was arbitrary, so $f \in C_h(X, \Ecal_L )$.
\qed
\enddemo

We prove that the map $h_L X \rightarrow \Xbar$ is a homeomorphism.

\proclaim{Proposition 2.9} Let $A$ and $B$ be subsets of a
proper metric space $(X,d)$, and suppose that there
is a constant $c>0$ such that $d(A \setminus B_r, B \setminus B_r) \geq cr$ for all
$r \geq 0$.  Then the function $\phi: X \rightarrow [0,1]$, defined by
$$ \phi(x) = \frac{d(x,A)}{d(x,A) + d(x,B)}, $$
is an element of $U(X, x_0)$, i.e. there is a $c_{\phi} > 0$ such that
$|\phi(x) - \phi(y)|\norm{x} \leq c_{\phi} d(x,y)$ for all $x,y \in X$.
\endproclaim

\demo{Proof} By the proof (not just the statement)
for the characterization of divergent systems with two members,
there is a number $C>0$ such that
$$d(x,A) + d(x,B) \geq \max{ \{ d(x,A), d(x,B) \} } \geq C\norm{x} $$
for all $x$ with $\norm{x} \geq 0$, that is for all $x \in X$.  So
$$|\phi(x) - \phi(y)|
\leq |\frac{d(x,A)}{d(x,A)+d(x,B)}-\frac{d(y,A)}{d(x,A)+d(x,B)}|
+
|\frac{d(y,A)}{d(x,A)+d(x,B)}-\frac{d(y,A)}{d(y,A)+d(y,B)}| $$
$$\leq  \frac{d(x,y)}{d(x,A)+d(x,B)}+\frac{d(y,A)}{d(y,A)+d(y,B)}
\frac{|d(y,A)-d(x,A)|+|d(y,B)-d(x,B)|}{d(x,A)+d(x,B)} $$
$$\leq  3 \frac{d(x,y)}{d(x,A) + d(x,B)} \leq \frac{3d(x,y)}{C\norm{x}}, $$
and the Proposition follows. \qed   \enddemo

\proclaim{Proposition 2.10} Let $(X,d)$ be a proper metric space, and suppose that
$A$ and $B$ are subsets of $X$.  Define $\nu X = \Xbar \setminus X$, and
set $A' = \overline{A} \cap \nu X$ and $B' = \overline{B} \cap \nu X$.  If there
exist $c,r_0 >0$ such that $d(A \setminus B_r, B \setminus B_r) \geq cr$ whenever
$r \geq r_0$, then $A' \cap B' = \varnothing$.  \endproclaim

\demo{Proof} Set $E = A \setminus B_{r_0}$ and $F = B \setminus B_{r_0}$. So
$d(E \setminus B_r, F \setminus B_r) = %
d((A \setminus B_{r_0}) \setminus B_r, (B \setminus B_{r_0}) \setminus B_r)$.
Then for $r \geq r_0$, we have
$d(E \setminus B_r, F \setminus B_r) = d(A\setminus B_r, B \setminus B_r) \geq cr$,
while if $r \leq r_0$, then
$ d(E \setminus B_r, F \setminus B_r) = d(A \setminus B_{r_0},B \setminus B_{r_0}) %
\geq cr_0 \geq cr$.
In any case, we have
$$ d(E \setminus B_r, F \setminus B_r) \geq cr \quad \text{for all \ } r \geq 0.$$

By the previous Proposition, the map $\phi(x) = \frac{d(x,E)}{d(x,E) + d(x,F)}$ lies
in $U(X, x_0)$; so $\phi$ extends to $\Xbar$, and we label this extension $\phi$ as
well.  Let $a \in A'$ and $b \in B'$.  So
$(a,b) \in \overline{A} \times \overline{B}$, and hence
there is a net $\{ (a_{\alpha}, b_{\alpha}) \}_{\alpha}$ of points
of $A \times B$ such that $a_{\alpha} \rightarrow a$ and $b_{\alpha} \rightarrow b$.
Since $a_{\alpha} \notin B_{r_0}$ and $b_{\alpha} \notin B_{r_0}$ eventually, and
hence $a_{\alpha} \in E$ and $b_{\alpha} \in F$ eventually, we have that
$\phi(a) - \phi(b) = \lim_{\alpha}( \phi(a_{\alpha}) - \phi(b_{\alpha}) ) = -1$.
That is, $\phi(a) \neq \phi(b)$, so $a \neq b$.  Since
$a \in A'$ and $b \in B'$ were arbitrary, have
$A' \cap B' = \varnothing$.\qed  \enddemo

\proclaim{Theorem 2.11} Let $(X,d)$ be a proper metric space.
Then the sublinear compactification $\Xbar$ is
homeomorphic to the Higson compactification $h_L X$ for the sublinear coarse
structure $\Ecal_L$ via a  homeomorphism extending the
identity on $X$.  \endproclaim
\demo{Proof}
We use $\nu_L X$ to denote the Higson corona associated with
the sublinear coarse structure; $\nu' X$ will denote the boundary of $\Xbar$, that is,
$\nu' X = \Xbar \setminus X$.

Let $\theta : h_LX \rightarrow \Xbar$ be this extension
of $\ident{X}$.  One can show that $\theta(\nu_L X) \subset \nu' X$.  Thus,
is suffices to show that the map is one-to-one on the corona.  So let
$x$ and $y$ be distinct points of $\nu_L X$.  So there are subsets $A$ and
$B$ of $X$ such that $x \in A'$, $y \in B'$, and
$\overline{A} \cap \overline{B} \cap \nu_L X = \varnothing$ (here the closure is
taken in $h_LX$).  This means that
$d(A \setminus B_r, B \setminus B_r) \geq cr$ eventually (for some c).  Thus, by
Proposition 2.10,
$\overline{A} \cap \overline{B} \cap \nu' X =\varnothing$
(closures take place in $\Xbar$).  But $\theta(x) \in \overline{A} \cap \nu' X$
and $\theta(y) \in \overline{B} \cap \nu' X$, so
$\theta(x) \neq \theta(y)$.\qed  \enddemo

REMARK. We have $C_h(X, \Ecal_L) = \overline{U(X,x_0)}$.  We will sometimes refer 
to a function $f \in C_h(X, \Ecal_L)$ as a linear Higson function on $X$.

\head \S3 The equality $\AN X=\dim\nu_LX$ \endhead

{\bf Asymptotic Assouad-Nagata dimension.}
We characterize the asymptotic Assouad-Nagata dimension using
a sequential formulation.

For a map $f: X \rightarrow Y$ between metric spaces, define
$$ \Lip{(f)} = \sup{ \{ \frac{d_Y (f(x),f(y))}{d_X (x,y)} : x, y \in X, x
\neq y \} }. $$
Note that this could be $\infty$.

For a map $p: X \rightarrow P$ to a uniform simplicial complex $P$ 
and a simplex $\Delta \subset P$, define
$$ D(p, \Delta) = \inf{ \{ \Lip{(g)} \mid %
g:p^{-1}(\Delta) \rightarrow \partial \Delta \text{ \ continuous}, %
g|_{p^{-1}\partial \Delta} = p|_{p^{-1}\partial \Delta} \} }. $$
It might be possible that there are no such $g$ with this property, or that
$\Lip(g) = \infty$ for all such $g$, in which case $D(p, \Delta) = \infty$.
Finally, we define
$$ D_n (p) = \sup{ \{ D(p, \Delta) : \Delta \subset P, \dim{\Delta} = n \} }.  $$
Based on the previous observations, this could be infinite as well.
In the event that $P$ is $n-$dimensional (the very case where we intend to
use these constructions), we will write $D(p)$ rather than $D_n (p)$.

The following lemma will be crucial for this section.

\proclaim{Lemma 3.1} For a metric space $(X,d)$, suppose that
$\AN{X} = n \geq 1$.
Then there is a sequence $\{ \lambda_m \}_{m=1}^{\infty}$
of positive numbers, a sequence
$\{ p_m : X \rightarrow P_m \}_m$ of maps to $n-$dimensional simlicial complexes $P_m$, 
and a number $C>0$ such that
\roster
\item $\la_m \rightarrow 0$ as $m \rightarrow \infty$;
\item $p_m$ is $\la_m-$Lipschitz and $C / \la_m -$cobounded;
\item $\ds \lim_{m \rightarrow \infty} \frac{D(p_m)}{\la_m} = \infty$.
\endroster \endproclaim
\demo{Proof} 
Set $a = (2n+3)^2$ for simplicity.  There is 
a $C > 0$, a positive integer $m_0$, and (for the sequence $\{ am \}_{m\geq m_0}$) 
a sequence of covers $\Ucal_m$ such that 
$L(\Ucal_m) > am$, $\mesh \Ucal_m \leq C a m$, and 
each ball $B_{am}(x)$ meets at most $n+1$ elements 
of $\Ucal_m$ for each $m \geq m_0$.     

Let $p_m : X \to \Nerve \Ucal_m$ be the projection to the nerve.  
It follows that $p_m$ is $\frac{1}{m}-$Lipschitz and 
$Cam-$cobounded [BD2].  Set $\la_m = \frac{1}{m}$ and 
$P_m = \Nerve \Ucal_m$.  

We now show that 
$ \limsup_{m \geq m_0} \frac{D(p_m)}{\la_m} = \infty$.  
To get a contradiction, suppose that there is a positive integer 
$m_1 \geq m_0$ and a  
$b > 1$ such that $\frac{D(p_m)}{\la_m} < b$ for all 
$m \geq m_1$.  Fix $m\geq m_1$, and let $\Delta$ be an $n-$simplex in 
$P_m$.  Thus, $D(p_m, \Delta) < b \la_m$, and so there is a 
$g_{\Delta} : p_m^{-1}(\Delta) \rightarrow \partial \Delta$ 
(depending on $m$) such that 
$\Lip (g_{\Delta}) < b\la_m$ and $g_{\Delta} = p_m$ on 
$p_m^{-1}(\partial \Delta)$.  
For $m \geq m_1$, we define a map 
$q_m: X \rightarrow (P_m)^{n-1}$ by 
$$ q_m(x) = \cases g_\Delta (x) & %
\text{if \ } x \in p_m^{-1}\Delta \text{\ for some $n-$dimensional simplex \ } \Delta \subset P_m \\ 
p_m (x) & \text{otherwise}. \endcases $$ 
Here $(P_m)^{n-1}$ denotes the $n-1$ skeleton of $P_m$.  
It is easy to see that $q_m$ is well-defined and 
$q_m = p_m$ on $p_m^{-1}( (P_m)^{n-1} )$.  

For $m \geq m_1$, define 
$\Vcal_m = \{ q_m^{-1} (\st v) : v \text{ is a vertex of } P_m \}$. 
We show $\mesh \Vcal_m \leq Cam$.  It suffices to show that 
$q_m^{-1} (\st v) \subset p_m^{-1} (\st v)$, where $v$ is a vertex of 
$(P_m)^{n-1}$ (the vertices of $(P_m)^{n-1}$ are the same of those of $P_m$, and 
can be identified with the elements of the cover $\Ucal_m$).  
Let $x \in q_m^{-1}(\st v)$, and we consider two cases.  First, 
suppose that $x \in p_m^{-1} ( (P_m)^{n-1} )$; then by definition, 
we have $p_m(x) = q_m (x) \in \st v$, and so 
$x \in p_m^{-1} (\st v)$.  Second, if 
$x \in p_m^{-1}(\interior \Delta)$ for some $n-$dimensional $\Delta$, then 
$q_m (x) \in \partial \Delta$ and since $q_m (x) \in \st v$, we have  
$v$ is a vertex of $\Delta$; as $p_m(x) \in \interior \Delta$, we have 
$p_m(x) \in \st v$, and so 
$x \in p_m^{-1} \st v$.  This proves the inclusion.  

We now show that 
$\Vcal_m$ is $\frac{m}{b(n+1)}-$Lipschitz.  Let $x \in X$.  
Let $U_0, U_1, \ldots, U_j$ be the distinct members of $\Ucal_m$ meeting 
$B(x, \frac{m}{b(n+1)})$.  By the construction of $\Ucal_m$, we must have $j \leq n$.  
Let $y \in B(x, \frac{m}{b(n+1)})$.  Note that 
$p_m (x)$ (and $p_m (y)$ as well) lies in a simplex 
whose vertices form a subset of $\{ U_i \}$.  
There is an $i$, $0 \leq i \leq j$, such that the $U_i-$th coordinate of $q_m (x)$ 
is at least $\frac{1}{n+1}$, i.e. $[q_m(x)]_{U_i} \geq \frac{1}{n+1}$.   
We consider two cases.  First, suppose that 
$p_m (x)$ and $p_m (y)$ lie in the $(n-1)-$skeleton of 
$P_m$; that is, $x,y \in p_m^{-1} ( (P_m)^{n-1} )$.  Second, if 
$p_m (x)$ or $p_m (y)$ lies in the interior of an $n-$simplex, 
then we must have $j \geq n$, and so $j=n$.  This means that 
$\{ U_i : i=0, 1, \ldots, n \}$ corresponds to an $n-$simplex $\Delta$ of $P_m$, 
and so $x,y \in p_m^{-1}(\Delta)$.  In either case, 
$$ |q_m (x) - q_m (y)|  < b\la_m d(x,y) \leq \frac{1}{n+1}, $$
and so $|[q_m(x)]_{U_i} - [q_m(y)]_{U_i}| < \frac{1}{n+1}$.   
Thus, $[q_m (y)]_{U_i} > 0$, or $y \in q_m^{-1} (\st U_i)$.  
This proves that $B(x, \frac{m}{b(n+1)}) \subset q_m^{-1}(\st U_i)$.  As $x \in X$ 
was arbitrary, we have $L(\Vcal_m) \geq \frac{m}{b(n+1)}$.  

So for $m \geq m_1$, $\Vcal_m$ is a cover of $X$, $\mult \Vcal_m \leq n $, 
$L(\Vcal_m) \geq \frac{1}{b(n+1)}m$, and 
$\mesh \Vcal_m \leq C a m $.  
Now let $r$ be a real number with 
$r \geq \frac{m_1}{b(n+1)}$, and choose an integer 
$m \geq m_1$ such that $m-1 \leq b(n+1)r \leq m$.  
Then $L(\Vcal_m) \geq \frac{m}{b(n+1)} \geq r$ and 
$\mesh \Vcal_m \leq C a m \leq C a(b(n+1)r + 1) \leq C a (b(n+1)r + m_1) \leq 2Cab(n+1)r$.  
Setting $r_0 = \frac{m_1}{b(n+1)}$ and $C_0 = 2Cab(n+1)$ gives 
$\AN X \leq n-1$, a contradiction.  

Thus, $\limsup \frac{D(p_m)}{\la_m} = \infty$.  Passing to a subsequence if 
necessary and relabeling, we have the desired result.  \qed \enddemo

\mk

{\bf Extensions of functions.}
Let $A\subset X$ be a closed subset. We call a neighborhood
$W\supset A$ {\it linear} if there is a constant $c>0$ such that
$d(A\setminus B_r,B\setminus B_r)\ge cr$ for all $r>0$ where
$B=X\setminus W$.

Suppose that $A_1 \subset A_2 \subset W$, where $A_1$ and
$A_2$ are closed subsets of $X$ and $W$ is a linear neighborhood
of $A_2$.  Then $W$ is a linear neighborhood of $A_1$ as well.
Let $A \subset W \subset X$, where $W$ is a linear neighborhood in $X$ of the closed
set $A$.  If $Y \subset X$ is a closed subset,
then $Y \cap W$ is a linear neighborhood in $Y$ of the closed set $Y \cap A$.

We shall extend the notation $U(X, x_0)$ as follows.  
For metric spaces $X$ and $Y$ (not necessarily proper), with $x_0 \in X$ and 
$\norm{\cdot} = d_X(\cdot, x_0)$, we say 
$f \in U(X, x_0, Y)$ if and only if 
$f$ is bounded, continuous, and there is a $c_f \geq 0$ such that 
$d_Y (f(x), f(y)) \norm{x} \leq c_f d_X (x,y)$ for all $x,y \in X$.  
In the event that $Y = \C$, we will omit $\C$ from the notation.  

REMARKS. Suppose that $f \in U(X,x_0, Y)$ and that 
$g:Y \to Z$ is a $\la-$Lipschitz map between metric spaces.  
Then 
$d_Z ((g \circ f)(x) , (g \circ f)(y))\norm{x} \leq \la c_f d_X(x,y)$, and 
so $g\circ f \in U(X, x_0, Z)$.  It is also clear that 
$f \in U(X, x_0, \R^{n+1})$ if and only if 
$f_i \in U(X,x_0,\R)$ for $1 \leq i \leq n+1$.

\proclaim{Proposition 3.2}
Let $g:X\to Y$ be an element of $U(X,x_0, Y)$ for proper metric spaces $X$ and 
$Y$.  
Then $W=g^{-1}(N_r(F))$ is a linear neighborhood
of $A=g^{-1}(F)$ for any closed subset $F\subset Y$ and $r > 0$.
\endproclaim
\demo{Proof} Let $c = c_g$.
Assume that $W$ is not linear. Then for each positive integer 
$n$, there is an $r_n > 0$ such that 
$d( (X \setminus W) \setminus B_{r_n}, A \setminus B_{r_n}) %
< \frac{1}{n}r_n$, and 
so there are $x_n,y_n\in X$ such that
$g(x_n)\in F$, $d(g(y_n),F) \geq r$, $\|x_n\|,\|y_n\|\ge r_n$, and
$d(y_n, x_n) < \frac{1}{n} r_n$.  
Then we obtain a contradiction:
$$
0<r\le |g(x_n)-g(y_n)|\le \frac{cd(x_n,y_n)}{\norm{x_n}} \leq %
c\frac{d(x_n, y_n)}{r_n} \le \frac{c}{n} \to 0.
$$ 
\qed \enddemo

\proclaim{Proposition 3.3} Let $q:X\to\R$ be a Higson function for the
sublinear coarse structure for a proper metric space $X$. Then
for every $\epsilon>0$ there is a linear neighborhood $W\supset
q^{-1}(0)$ such that $W\subset q^{-1}(-\epsilon,\epsilon)$.
\endproclaim

\demo{Proof} Since $q \in C(\Xbar)$, there is a  function
$g \in U(X,x_0)$ with $|g-q|<\epsilon/4$. Take $W=g^{-1}(-\epsilon/2,\epsilon/2)=
g^{-1}(N_{\epsilon/4}([-\epsilon/4,\epsilon/4]))$.
By Proposition 3.2 it is a linear neighborhood of
$g^{-1}([-\epsilon/4,\epsilon/4] \supset q^{-1}(0)$.\qed  \enddemo

\proclaim{Proposition 3.4} Let $(X,d)$ be a proper metric space with
basepoint $x_0$.  Let $A$ be a closed subset of $X$ containing
the basepoint, and let $W$ be an open linear neighborhood of $A$.  
Suppose that $f \in U(W,x_0)$ and let
$\bar f$ be an extension of $f|_A$ which is a linear Higson
function on $X$. Then for every $\epsilon>0$ there is a
$g \in U(X,x_0)$  which is $\epsilon$-close to $\bar f$ and extends $f|_A$.  
If $f$ and $\bar f$ are real-valued, then $g$ can be take to be real-valued 
as well.  
\endproclaim
\demo{Proof} We first construct a map $\hat f \in U(X, x_0)$
with $\hat{f}|_A = f|_A$.  Set $B = X \setminus W$.  Let $\phi:X \rightarrow \R$
be defined by $\phi(x) = \frac{d(x,B)}{d(x,A)+d(x,B)}$.  Since $W$ is a linear neighborhood,
there is a $c$ such that
$d(A \setminus B_r, B \setminus B_r) \geq cr$ for all $r \geq 0$; thus,
$\phi \in U(X,x_0)$ by Proposition 2.9.  Define $\hat{f}: X \rightarrow \C$ by taking
$\hat{f}(x) = \phi(x)f(x)$ if $x \in W$, and setting $\hat{f}(x) = 0$ otherwise.
Clearly $\hat f$ is bounded and $\hat{f}|_A = f|_A$.

Let $x,y \in X$.  We consider four cases.  First, if
$x, y \in W$, then
$| \hat{f}(x) - \hat{f}(y)| \cdot \norm{x} = | \phi(x)f(x) - \phi(y)f(y) | \norm{x} %
\leq (c_{f} + \norm{f} c_{\phi})d(x,y)$,
where $c_{f}$ and $c_{\phi}$ are the appropriate constants for $f$ and $\phi$.
If $x \in W$ and $y \notin W$, then
$|\hat{f}(x) - \hat{f}(y)|\norm{x} = |\phi(x)f(x)| \norm{x} %
\leq \norm{f}|\phi(x)-\phi(y)| \norm{x} \leq \norm{f}c_{\phi} d(x,y).$
If $x \notin W$ and $y \in W$, then a similar argument shows that
$|\hat{f}(x) - \hat{f}(y)|\norm{x} \leq \norm{f}c_{\phi} d(x,y)$.  
Finally, $|\hat{f}(x) - \hat{f}(y)| \norm{x} = 0$ when $x,y \notin W$.
Thus, we have that
$|\hat{f}(x) - \hat{f}(y)| \norm{x} \leq ( c_f + \norm{f}c_{\phi} )\norm{x} $ for all
$x$ and $y$ in $X$.  By the remark before Proposition 2.7, $\hat{f}$ is continuous
everywhere except possibly at $x_0$; but $W$ is an open neighborhood
containing $x_0$, so $\hat{f}$ is continuous at $x_0$ since
$\hat{f} = \phi f$ on $W$.  Hence $\hat{f} \in U(X,x_0)$.

Let $\tilde f \in U(X,x_0)$ be an $\epsilon/2$-approximation of
$\bar f$.  We consider the function
$q=\hat f - \bar{f} :X \to \C$.  By Proposition 3.3 there is a
linear neighborhood $W_0$ of $A\subset q^{-1}(0)$ such that
$W_0\subset q^{-1}(-\epsilon/2,\epsilon/2)$. Let $\psi_1$, $\psi_2$ be defined
by
$$ \psi_1 (x)= \frac{d(x,X \setminus W_0)}{d(x,A)+d(x,X\setminus W_0)} \quad \text{and} \quad %
\psi_2 (x)= \frac{d(x,A)}{d(x,A) + d(x, X\setminus W_0)};  $$
so $\psi_1, \psi_2 \in U(X,x_0)$ by Proposition 2.9.  
We define $g=\psi_1 \hat f + \psi_2 \tilde f$. Then $g \in U(X,x_0)$ since
$U(X,x_0)$ is an algebra. Note that
$$
|g-\bar f|=|\psi_1 \hat{f} + \psi_2 \tilde f- (\psi_1+\psi_2)\bar f | %
\le |\psi_1 (\hat f - \bar f )| + |\psi_2 ||\tilde f-\bar f|\le\epsilon/2+\epsilon/2.
$$
Finally, it is clear that $g|_A = \hat{f}|_A = f|_A$.
\qed
\enddemo

\proclaim{Proposition 3.5}  Suppose that $u,v \in U(X,x_0)$ are nonnegative,
real-valued functions, and
$u(x)+v(x) \geq \delta >0$ for all $x \in X$.  Then
$\frac{u}{u+v} \in U(X,x_0)$ as well.
\endproclaim
\demo{Proof}
$$ |\frac{u(x)}{u(x)+v(x)} - \frac{u(y)}{u(y)+v(y)} | \norm{x}
 \leq \left( \frac{|u(x)-u(y)|}{u(x)+v(x)} +
|u(y)| |\frac{1}{u(x)+v(x)}-\frac{1}{u(y)+v(y)} | \right) \norm{x}
$$
$$\leq \frac{|u(x)-u(y)|}{u(x)+v(x)} \norm{x}
 + u(y) \frac{|u(x)-u(y)|+|v(x)-v(y)|}{(u(x)+v(x))(u(y)+v(y))} \norm{x}
\leq \frac{2c_u + c_v}{\delta} d(x,y).  $$
\qed \enddemo

Proposition 3.5 also holds if we replace $U(X,x_0)$ by $U(X,x_0,\R)$ throughout 
the statement.  In what follows, $f_i$ will indicate the $i$th component of a
function $f:X \rightarrow \R^{n+1}$.  

\proclaim{Lemma 3.6} Let $(X,d)$ be a proper metric space with
basepoint $x_0$.  Let $A$ be a closed subset of $X$ containing
the basepoint, and let $W$ be an open linear neighborhood of $A$.
Let $\Delta$ denote the standard $n-$simplex.
Suppose that $f \in U(W, x_0, \partial \Delta)$
and $g:X\rightarrow \partial \Delta$ is a
continuous extension of $f|_A$ such that
each component of $g$ is a linear Higson function.
Then there is an 
$h \in U(X,x_0, \partial \Delta)$ which extends $f|_A$. \endproclaim

\demo{Proof} Looking at components, by Proposition 3.4 we have
that there are $q_i \in U(X,x_0, \R)$ such that $q_i |_A = f_i |_A$ and 
$\norm{q_i - g_i} \leq \frac{1}{3(n+1)}$.  
Thus,
$$ |1 - \sum_i |q_i| | = |\sum_i |g_i| - \sum_i |q_i| |
\leq \sum_i | |g_i| - |q_i| | \leq \frac{1}{3}, $$
and so $\sum_i |q_i(x)| \geq 2 / 3$ for all $x \in X$.  Define
$q': X \rightarrow \Delta$ by
$q'(x) = \left( \frac{|q_i (x)|}{\sum_j |q_j (x) |} \right)_{i=1}^{n+1}. $
This map is well-defined.  Also, since $q_i \in U(X,x_0,\R)$, we have
$|q_i| \in U(X, x_0,\R)$, and hence $q_i' \in U(X,x_0,\R)$ by Proposition 3.5.

Now, fix $x \in X$; note that there is a $j$
such that $g_j(x) = 0$.  So
$|q_j (x)| \leq \norm{q_j - g_j} \leq \frac{1}{3(n+1)}$.  Thus,
$\frac{|q_j(x)|}{\sum_k |q_k(x)|} \leq \frac{1}{2(n+1)}$, and so
$$d(q'(x), b) =
\sqrt{\sum_{1\leq i \leq n+1} | \frac{|q_i(x)|}{\sum_k |q_k(x)|} - \frac{1}{n+1} |^2} %
\geq | \frac{|q_j(x)|}{\sum_k |q_k(x)|} - \frac{1}{n+1}  | \geq \frac{1}{2(n+1)}, $$ 
where $b$ is the barycenter of $\Delta$.  
That is, $q'$ maps $X$ into $\Delta \setminus B_{\frac{1}{2(n+1)}}(b)$.  Let
$r: \Delta \setminus B_{\frac{1}{2(n+1)}}(b) \to \partial \Delta$ be a
$\la-$Lipschitz retract onto $\partial \Delta$.  Finally, define
$h:X \to \partial \Delta$ by $h = r \circ q'$.  So
$h|_A = r \circ q' |_A = r \circ f|_A = f|_A$ since $r$ is a retract.  Also, since
$q_i' \in U(X, x_0, \R)$ for all $i$, we have
$q' \in U(X,x_0, \Delta \setminus B_{\frac{1}{2(n+1)}}(b) )$, and so
$h \in U(X,x_0, \partial \Delta )$ by the remarks preceding Proposition 3.2. \qed
\enddemo

\mk

{\bf The inequality $\AN X\le\dim\nu_LX$.}

\mk

We recall that a metric space $(X,d)$ is called {\it cocompact}
if there is a compact subset $K$ of $X$ such that
$X = \cup_{\g \in \operatorname{Isom}(X)} \g (K)$, where
$\operatorname{Isom}(X)$ is the set of all isometries of $X$.

\proclaim{Theorem 3.7} Let $X$ be a cocompact, connected, proper metric space which has
finite asymptotic Assouad-Nagata dimension.
Then  $\dim \nu_L X \geq \AN X$.  \endproclaim
\demo{Proof} Set $n = \AN X$.  If $n = 0$, the inequality is immediate.  We shall henceforth
assume that $n > 0$, and so in particular $X$ is not compact.
To get a
contradiction, assume that $\dim \nu X \leq n-1$.  By hypothesis,
there is a compact subset $K$ of $X$ such that
$X = \cup_{\g \in \G} \g (K)$, where $\G = \operatorname{Isom}(X)$.

Let  $\{ \lambda_m \}_{m=1}^{\infty}$ be a sequence
of positive numbers, let
$\{ p_m : X \rightarrow P_m \}_m$ be a sequence of maps to
$n-$dimensional polyhedra, and let $C>0$ be a constant such that
(1) - (3) of Lemma 3.1 hold.  Without loss of generality, one may take
$C > \diam K$.  Also, passing to subsequences if necessary,
we may assume that $\la_1 \leq 1$ and 
$\la_{i+1} \leq \la_i / 25$.  For every $i$, we take
$\Delta_i \subset P_i$ with $D(p_i, \Delta_i) \geq D(p_i)/2$.  Let
$h_i : \Delta_i \to \Delta$ be an isometry to the standard
$n-$simplex.

Let $\g_i \in \G$ be an element with
$p_i^{-1}(\Delta_i) \cap \g_i K \neq \varnothing$; let
$x_i \in p_i^{-1}(\Delta_i) \cap \g_i K$.  Choose $y_i$ with
$\norm{y_i} = \frac{3C}{\la_i}$, and let
$\alpha_i$ be such that $y_i \in \alpha_i K$.  Define
$$ A_i = \alpha_i \g_i^{-1} p_i^{-1}(\Delta_i) \quad \text{and} \quad
B_i = \alpha_i \g_i^{-1} p_i^{-1}(\partial \Delta_i) $$
for $i=1,2,\ldots$.  We also define $A_0 = B_0 = \{ x_0 \}$.
Note that $\diam A_i \leq \frac{C}{\la_i}$ when $i \geq 1$.  

For $a \in A_i$,
$$d(a, y_i) \leq d(a, \alpha_i \g_i^{-1} x_i) + d(\alpha_i \g_i^{-1} x_i, y_i)
\leq \frac{C}{\la_i} + C \leq \frac{2C}{\la_i}.  $$
So for $a \in A_i$, we have
$\frac{3C}{\la_i} = \norm{y_i} \leq \norm{a} + d(a, y_i) \leq \norm{a} + \frac{2C}{\la_i}$
and so $d(x_0,A_i) \geq \frac{C}{\la_i}$.  Also, for $a \in A_i$,
$\norm{a} \leq \norm{y_i} + d(y_i, a) \leq \frac{5C}{\la_i}, $
and so $\norm{A_i} := \sup_{a \in A_i} \norm{a} \leq \frac{5C}{\la_i}$ for $i \geq 1$.
Note that the $A_i$ ( $i \geq 0$) are disjoint since
$$\norm{A_i} \leq \frac{5C}{\la_i} \leq \frac{C}{5\la_{i+1}} <
\frac{C}{\la_{i+1}} \leq d(x_0,A_{i+1}) $$
for all $i \geq 1$ and since $x_0 \notin A_i$ for all $i\geq 1$ (because $d(x_0, A_i) > 0$).  Set
$A = \coprod_{i\geq 0} A_i$ and $B = \coprod_{i\geq 0} B_i$.  
$A$ and $B$ are closed in $X$.

Let $q_i = h_i p_i \g_i \alpha_i^{-1} |_{A_i}$ for $i \geq 1$ and define
$q_0 : A_0 \rightarrow \Delta$ by $q_0 (x_0) = (1,0,0, \ldots, 0)$ (the image does not matter,
however).  Now take $Q = \coprod_{i\geq 0} q_i : A \rightarrow \Delta$.
So
$\norm{Q(x) - Q(y)} \norm{x} \leq (\Lip{q_i}) \norm{A_i} d(x,y) \leq \la_i \frac{5C}{\la_i} d(x,y)$
when $x,y \in A_i$ for $i \geq 1$.  For $x,y \in A_0 = \{ x_0 \}$, we just
have $\norm{Q(x) - Q(y)} \norm{x} = 0$.

Now suppose that $x \in A_j$ and $y \in A_i$, where $j \neq i$ and $i,j \geq 1$.
We first prove this in the case that $j>i$.  So
$\norm{ Q(x) - Q(y) } \norm{x} \leq \sqrt{n+1} \norm{A_j} \leq \frac{5\sqrt{n+1}C}{\la_j} $
and since $\la_j \leq \la_i / 25$, we have
$d(x,y) \geq d(x_0, A_j) - \norm{A_i} \geq \frac{C}{\la_j} - \frac{5C}{\la_i}
 \geq \frac{C}{\la_j} \left( 1 - \frac{1}{5} \right) = \frac{4C}{5\la_j}.$

Thus, we have
$$ \norm{Q(x) - Q(y) } \norm{x} \leq 5 \sqrt{n+1} (\frac{5}{4}) d(x,y) = \frac{25}{4}\sqrt{n+1}d(x,y). $$
If $j < i$, then $\norm{x} \leq \norm{y}$, and so
$ \norm{Q(x) - Q(y) } \norm{x} \leq \norm{Q(y)-Q(x)} \norm{y}$,
and thus by appealing to the case just considered, we have
$\norm{Q(x)-Q(y)} \leq \frac{25}{4}\sqrt{n+1}d(y,x)$.  
Now suppose $x \in A_j$ and $y \in A_i$, where $i$ and $j$ are distinct and
either $j=0$ or $i=0$.  So either $x = x_0$ or $y = x_0$, and in either case it
is easy to check that $ \norm{Q(x)- Q(y)} \norm{x} \leq \sqrt{n+1} d(x,y)$.

Thus, for
$c_Q = \max{ \{ 5C, \frac{25}{4}\sqrt{n+1}  \} }$, we have that
$\norm{Q(x)-Q(y)} \norm{x} \leq c_Q d(x,y)$ for all $x,y \in A$.
That is, $Q \in U(X,x_0, \Delta)$.

We now construct a map $f:W \to \partial\Delta$, where 
$W$ is a linear neighborhood of $B$ in $A$, 
$f$ extends $Q|_B$, and $f \in U(W, x_0, \partial \Delta)$.  
Set $\e = \frac{1}{2(n+1)}$ and set $W = Q^{-1}(N_{\e}(\partial\Delta))$.  
$W$ is a linear neighborhood of $B$ by Proposition 3.2, 
and $Q|_{W} \in U(W, x_0, \Delta)$.   Let  
$b \in \Delta$ be the barycenter and 
note that $B_{\e}(b) \cap N_{\e}(\partial\Delta) = \varnothing$ (neighborhood is taken in 
$\Delta$).  There is a 
Lipschitz retraction $r : \Delta \setminus B_{\e}(b) \to \partial \Delta$, and 
so $f:= r \circ Q|_W : W \to \partial\Delta$ is in $U(W,x_0,\partial\Delta)$ by 
the remarks preceding Proposition 3.2.  Clearly $f|_B = Q|_B$.    

Also $Q$, when restricted to $\coprod_{i \geq 0} B_i$, is a map into $\partial \Delta$.
Each $Q_i$ extends to a map $\overline{B}$ into $\R$, and so
$Q$ extends to a map $G'' : \overline{B} \to \R^{n+1}$.  Since $G''$ is continuous,
$G''(\overline{B}) \subset \overline{G''(B)} \subset \overline{ \partial\Delta} = \partial\Delta$,
and so we view $G''$ as a map from $\overline{B}$ to $\partial\Delta$.
Since $\nu B$ is closed in $\nu A $, where
$\dim \nu A \leq \dim \nu X \leq n-1$, we have an extension of $G''|_{\nu B}$
to $\nu A$; thus, we have an extension $G':\overline{B} \cup \nu A \to \partial\Delta$
of $G''$.  Since $\partial \Delta$ is an
absolute neighborhood retract, there is a continuous extension $G:V \to \partial\Delta$
of $G'$ to a neighborhood $V$ of
$\overline{B} \cup \nu A$ in $\overline{A}$.  Thus, there is an
$m_0$ such that $A' := A_0 \coprod (\coprod_{i \geq m_0} A_i)$ is a subset of
$V$; set $g= G|_{A'}: A' \rightarrow \partial\Delta$ and note that
$g|_{B'} = Q|_{B'}$, where
$B' = B_0 \coprod (\coprod_{i \geq m_0} B_i)$.  Also, $A'$ is closed in $A$ and so
$\overline{A'} \subset V$.  So $g$ extends to a continuous
map on $\overline{A'} \subset V$, and
since the sublinear Higson compactification of $A'$ is homeomorphic to
the closure of $A'$ in $\overline A$, we have that each component of
$g$ is a Higson function on $A'$.

We now restrict our attention to $A'$.  Note that
$W \cap A'$ is a linear neighborhood of $B \cap A' = B'$ in $A'$.
We have that $f|_{A'\cap W}: A' \cap W \to \partial\Delta$ is an element of
$U(A' \cap W , x_0, \partial \Delta)$, each component of $g$ is a Higson function,
and $g|_{B'} = Q|_{B'} = f|_{B'}$.
By lemma 3.6, we have that there is an
$h: A' \rightarrow \partial\Delta$ such that $h$ extends $f|_{B'}$ and for which
there is a $c_h$ such that
$\norm{h(x) - h(y)} \norm{x} \leq c_h d(x,y) $
for all $x,y \in A'$.

We now look at $A_i$ for $i \geq m_0$.
So $h|_{A_i} : A_i \rightarrow \partial\Delta$ extends
$f|_{B_i} = h_i p_i \gamma_i \alpha_i^{-1}|_{B_i}$ and
$$ \frac{C}{\la_i} \norm{h(x) - h(y)} \leq \norm{h(x) - h(y)} \norm{x} \leq c_h d(x,y) $$
whenenver $x,y \in A_i$.  Thus, $h|_{A_i}$ is $\frac{c_h \la_i}{C} -$Lipschitz.
So
$$ \frac{D(p_i)}{2} \leq D(p_i, \Delta_i) \leq \Lip(h_i^{-1}h|_{A_i}\alpha_i \gamma_i^{-1}) %
= \Lip( h|_{A_i} )  \leq \frac{c_h \la_i}{C} $$
and thus
$\frac{2 c_h}{C} \geq \frac{D(p_i)}{\la_i} \to \infty$,
a contradiction.
\qed
\enddemo

When we apply Theorem 3.7 to the Cayley graph of a finitely generated
group we obtain the following.

\proclaim{Corollary 3.8} For a finitely generated group $\G$ with word metric,
$\dim \nu_L \G \geq \AN{\G}$ provided $\AN{\G} < \infty$.  \endproclaim

EXAMPLE. Consider the parabolic region
$X = \{ (x,y) \in \R^2 : x \geq 0, |y| \leq \sqrt{x} \}$, which we
will equip with the (restricted) Euclidean metric.
Let $i: [0, \infty) \to X$ be the map
$i(x) = (x, 0)$.  Taking the usual metric on $[0, \infty)$,
it is not hard to show that $i$ is a coarse equivalence
for the sublinear coarse structures, and hence there is
a homeomorphism $\nu_L [0, \infty) \to \nu_L X$.  In particular,
$\dim \nu_L X = 1$.  But $\AN X = 2$.  
Since $X$ is a connected proper metric space, this shows that
we can not drop the requirement in the Theorem that the space be cocompact.

\mk

{\bf  The inequality $\AN X\ge\dim\nu_LX$.}

\mk

If $\Ucal$ is a cover of $X$ and $A \subset X$, then we write
$\Ucal_A = \{ U \in \Ucal : U \cap A \neq \varnothing \}$.

\proclaim{Lemma 3.9} Let $\Ucal$ and $\Vcal$ be covers of $X$, and suppose that $\Ucal$ refines $\Vcal$.  Let
$K$ be a subset of $X$.  So for each $U \in \Ucal$ with $U \cap (X \setminus K) \neq \varnothing$,
there is a $V_{U} \in \Vcal$ with $U \subset V_U$.  For $V \in  \Vcal$, set 
$$V' = [V \cap (X \setminus K)] \bigcup %
[\bigcup_{ U \in \Ucal_{X \setminus K}, V_U = V } U ] $$
and define 
$$ \Wcal = \{ U \in \Ucal | U \subset K \} \cup \{ V' | V \in \Vcal_{X \setminus K} \}. $$
Then
\roster
\item $\Wcal$ is a cover of $X$;
\item $\mult\Wcal \leq n+1$ if $\mult \Ucal \leq n+1$ and $\mult \Vcal \leq n+1$;
\item $\Wcal$ refines $\Vcal$;
\item $\Ucal$ refines $\Wcal$;
\item if $W \in \Wcal$ and $W \subset K$, then $W \in \Ucal$;
\item if $\Ucal$ and $\Vcal$ are open covers and $K$ is closed in $X$, then $\Wcal$ is also an open cover;
\item if $V \in \Vcal$ and $V \cap K = \varnothing$, then $V = V' \in \Wcal$.
\endroster
\endproclaim
The proof is strightforward.  For the above $\Wcal$ we will write
$$ \Wcal = \Ucal *_K \Vcal .$$
The following Theorem is a modification of Lemma 2.9 of \cite{DKU}.

\proclaim{Theorem 3.10}  Let $(X,d)$ be a proper
metric space.  Then
$\dim \nu_L X \leq \AN X$.  \endproclaim

\demo{Proof}
We write $\nu X = \nu_L X$ and use $\Bd(x,r)$ to denote a closed ball of radius $r$.  
Set $n = \AN X$ (if $\AN X = \infty$, the inequality is immediate).
As $\{ \widetilde{U} \cap \nu X : U \subset X \text{\ open} \} $
is a basis for $\nu X$, it suffices to prove that each cover of the form
$\{ \widetilde{U}_i \cap \nu X : 1\leq i \leq m \}$, where each $U_i \subset X$ 
is open, 
admits a finite refinement of multiplicity $\leq n+1$.

So let $\{ \widetilde{U}_i \cap \nu X : 1\leq i \leq m \}$
be a cover of $X$, and set $\Ucal = \{ U_i : 1 \leq i \leq m \}$.
Since $\AN X =n$, there exist $C > 0$ and $r_{-1} > 0$ such that whenever $r \geq r_{-1}$,
there is an open cover $\Ucal(r)$ of $X$ satisfying
$\mult \Ucal(r) \leq n+1$, $\mesh \Ucal(r) < Cr$, and $L(\Ucal(r)) > r$.  Without loss
of generality, we take $C > 1$.
Also, there is a $D>0$ and an $r_{-2} > 0$ such that
$L^{\Ucal}(x) \geq D \norm{x}$ whenever $x$ is such that
$\norm{x} \geq r_{-2}$.  We may take $D< 1$.

Now, choose $r_0 > \max{ \{ r_{-2}, r_{-1}, 1 \} }$.  Define
$r_i = (\frac{C}{D})^i r_0$ for $i \geq 1$.  Observe that
$r_{i+1} = \frac{C}{D} r_i > r_i > r_0 > 1$.  
Since $r_i > r_0 > r_{-1}$ for $i \geq 1$, there is a  cover
$\Ucal_i$ of $X$ such that $\mult \Ucal_i \leq n+1$, $\mesh \Ucal_i < Cr_i$,
and $L(\Ucal_i) > r_i$.

Define $\Vcal_1 = \Ucal_1$, and note that $\mesh \Vcal_1 < C r_1$.  Now,
supposing we have defined $\Vcal_1, \Vcal_2, \ldots, \Vcal_i$ satisfying
$\mesh \Vcal_j < Cr_j < r_{j+1}$ for all $1 \leq j \leq i$, then
$\Vcal_i$ refines $\Ucal_{i+1}$, and so we can define
$$ \Vcal_{i+1} = \Vcal_i *_{\Bd(x_0, 2r_{i+2})} \Ucal_{i+1}. $$
By Lemma 3.9, $\Vcal_{i+1}$ refines $\Ucal_{i+1}$, and so
$\mesh \Vcal_{i+1} < Cr_{i+1}$.  Thus, we have constructed $\Vcal_i$ for
all positive integers $i$.

Set $\Vcal = \liminf_i \Vcal_i = \cup_s \cap_{t \geq s} \Vcal_t$.
We now investigate some properties of $\Vcal$.

Using the definition of $\Vcal_{i+1}$, it is easy to show 
that if $U \in \Vcal_i$ and $U \cap \Bd(x_0, r_{i+2}) \neq \varnothing$,
then $U \in \Vcal_{i+1}$.  We conclude that
$\{ U \in \Vcal_i : U \cap \Bd(x_0,r_{i+2}) \neq \varnothing \} \subset \Vcal$.  
As $\Vcal_i$ is a cover of $X$ and hence of $\Bd_{r_{i+2}}$, we have that
$\Vcal$ covers $\Bd_{r_{i+2}}$;as $i$ here is arbitrary, $\Vcal$ covers $X$.

We now show that if $V \in \Vcal$ and $V \cap \Bd_{r_{i+1}} \neq \varnothing$,
then $V \in \Vcal_{i-1}$.  First suppose that $V \in \Vcal_i$ and
$V \cap \Bd_{r_{i+1}} \neq \varnothing$; then $\mesh \Vcal_i < Cr_i < r_{i+1}$ 
implies that $V \subset \Bd_{2r_{i+1}}$ 
and so $V \in \Vcal_{i-1}$ by (5) of the lemma.  Now suppose 
$V \in \Vcal$.  This means there is an
$s \geq i-1$ such that $V \in \Vcal_s$.  Applying the result we just found and
proceeding inductively, one can show
that $V \in \Vcal_j$ for all $j$ such that $i-1 \leq j \leq s$.

We show that $\Vcal_i$ refines $\Vcal$ for all $i \geq 1$.  We
know by (4) of the lemma that $\Vcal_i$ refines $\Vcal_{i+1}$
for all $i \geq 1$.  Fixing $i$, let $V \in \Vcal_i$.  Choose
$j \geq i$ such that $V \cap \Bd_{r_{j+2}} \neq \varnothing$.  As
$\Vcal_i$ refines $\Vcal_j$, there is a $U \in \Vcal_j$ such that
$V \subset U$.  Also,
$U \cap \Bd_{r_{j+2}} \supset V \cap \Bd_{r_{j+2}} \neq \varnothing$.  Thus,
$V \subset U \in \Vcal$.  So $\Vcal_i$ refines $\Vcal$.

Since each $\Vcal_i$ has multiplicity $\leq n+1$ for each $i$, it is clear 
from the definition that $\Vcal$ has multiplicity $\leq n+1$.  

Set $\Wcal = \Vcal_{X \setminus \Bd_{r_2}}$.
We have that $(\Vcal_i)_{X \setminus \Bd_{r_2}}$ refines $\Wcal$.

We show that $\Wcal$ refines $\Ucal$.  Let $W \in \Wcal$.  So
there is an $x \in W$ such that $\norm{x} > r_2$.  Take
$i = \max{ \{ j : \norm{x} > r_j \} }$, and note that
$i \geq 2$.  Thus, $\norm{x} \leq  r_{i+1}$, or $x \in \Bd_{r_{i+1}}$.  Hence
$W \cap \Bd_{r_{i+1}} \neq \varnothing$.  Since $W \in \Vcal$, we have
$W \in \Vcal_{i-1}$.  So
$$ \diam W < Cr_{i-1} = D (\frac{C}{D}) r_{i-1} = Dr_i \leq D\norm{x} \leq L^{\Ucal}(x), $$
and so there is a $U \in \Ucal$ with $W \subset U$.

We show that $L^{\Wcal}: X \to [0, \infty)$ is at least linear.  Set
$a = 3r_2$, and let $x$ be an element of $X$ with $\norm{x} \geq a = 3r_2$.
Set $i = \max{ \{ j : 3r_{j+1} \leq \norm{x} \} }$, and note that
$i \geq 1$ and $3r_{i+1} \leq \norm{x} < 3r_{i+2}$.  As
$L(\Ucal_i) > r_i$, there is a $U \in \Ucal_i$ such that
$B(x,r_i) \subset U$.  Since $\norm{x} \geq 3r_{i+1}$ and
$\diam U \leq \mesh \Ucal_i < Cr_i < r_{i+1}$, we have that
$U \subset X \setminus \Bd(x_0, 2r_{i+1})$.  By definition of
$\Vcal_i$, and by (7) of the lemma, we have that $U \in \Vcal_i$. In fact,
as $\norm{x} > r_2$, we have $U \in (\Vcal_i )_{X \setminus \Bd_{r_2}}$.
Since $(\Vcal_i )_{X \setminus \Bd_{r_2}}$ refines $\Wcal$, we have that
$$ L^{\Wcal}(x) \geq d(x,X \setminus U) \geq r_i.  $$
But $\norm{x} < 3 r_{i+2} = 3\frac{C^2}{D^2} r_i, $
so $L^{\Wcal}(x) > \frac{D^2}{3C^{2}} \norm{x}$.  Therefore,
$L^{\Wcal}$ is at least linear.

To summarize, $\Wcal$ covers $X \setminus \Bd_{r_2}$, $\mult \Wcal \leq n+1$,
$\Wcal$ refines $\Ucal$, and $L^{\Wcal}$ is at least linear.  Thus,
for $W \in \Wcal$, there is a $U_W \in \Ucal$ for which
$U \subset U_W$.  So for each $1 \leq i \leq m$, we define
$ W_i = \cup_{U_W = U_i} W$.

Now set $\Wcal' = \{ W_i \}$.  It follows that $\Wcal$
refines $\Wcal'$ and $\Wcal'$ has multiplicity $\leq n+1$.  Thus, $L^{\Wcal'} \geq L^{\Wcal}$
and hence $L^{\Wcal'}$ is at least linear.  As a consequence, if
we define $\widetilde{\Wcal}' = \{ \widetilde{W}_i \cap \nu X \}$, we have
that $\widetilde{\Wcal}'$ is a cover of $\nu X$.  Since
$\Wcal'$ refines $\Ucal$, we have that $\widetilde{\Wcal}'$ refines
$\widetilde{\Ucal}$.  Finally, as $\Wcal'$ has multiplicity $\leq n+1$, so
does $\widetilde{\Wcal}'$.\qed  \enddemo

\proclaim{Corollary 3.11} For a cocompact connected proper metric space, 
$\AN{X}=\dim\nu_LX$ provided $\AN{X}<\infty$.
\endproclaim

In particular, we conclude that $\AN \G = \dim\nu_L \G$ for a finitely generated 
group $\G$ with $\AN \G < \infty$.

\head \S4 Morita type theorem for Assouad-Nagata dimension\endhead

Let $K\subset S^n$ be a compact set in the unit sphere $S^n\subset\R^{n+1}$.
The open cone
$OK$ by definition is the union of rays through $K$ issuing from the origin
with the metric restricted from $\R^{n+1}$. The open cone admits a natural
compactification by $K$, which we will denote by $\overline{OK}$.

\proclaim{Proposition 4.1} The natural compactification of the open cone $OK$
is dominated by the sublinear compactification.
\endproclaim
\demo{Proof} It suffices to show that if two sets $A,B\subset OK$ do not
intersect at the cone boundary then they are divergent in the sublinear coarse
structure, but this is obvious.\qed
\enddemo 

\proclaim{Lemma 4.2} Let $X$ be a proper metric space. Then there is an embedding
$$ \nu_LX\times[0,1]\to\nu_L(X\times\R_+).$$
\endproclaim
\demo{Proof}
Let $d_0:X\to\R_+$ be the distance function to the base point $x_0\in X$, i.e.,
$d_0(x)=\|x\|$.
Then the map $d_0\times 1:X\times\R_+\to\R_+\times\R_+$ is a coarse morphism
for the sublinear coarse structures and hence is extendible to the sublinear
Higson compactifications:
$$\gamma=\overline{(d_0\times 1)}: h_L(X\times\R_+)\to h_L(\R_+\times\R_+).$$

Let $K$ be the arc in $S^1\subset\R^2$ from 0 to $\pi/4$.
In view of Proposition 4.1 we have a 
natural map $\phi':h_L OK\to \overline{OK}$. 
Also, $\tan$ can be defined on $\overline{OK} \setminus \{0 \}$. 
Let $W=(d_0\times 1)^{-1}(OK)$. 
Define $\phi:h_L W \setminus \{ (x_0, 0) \} \to [0,1]$
as the composition $(\tan)\circ\phi'\circ\gamma$ restricted to 
$h_L W \setminus \{ (x_0, 0) \}$.  Note that, 
by Corollary 2.6, $\nu_LW\subset\nu_L(X\times\R_+)$.

The restriction of the projection $X\times\R_+\to X$ to $W$ is a coarse morphism
and hence it
defines a map $\psi:\nu_LW\to\nu_LX$. We show that the map
$\Phi=(\phi|_{\nu_L W},\psi):\nu_LW\to [0,1]\times\nu_LX$ is a homeomorphism.

First we note that for every $t\in K$ the preimage
$$X_t=(d_0\times 1)^{-1}(Ot)=\{(x,\|x\|\tan{(t)})\mid x\in X\}\subset X\times\R_+$$
is coarsely equivalent to $X$. This implies that the map $\Phi$
takes $\nu_LX_t$ onto $\tan(t)\times\nu_LX$ homeomorphically. Thus,
$\Phi$ is onto. 

It remains to show that
$\Phi^{-1}(\tan(t)\times\nu_LX)=\nu_LX_t$, or equivalently
$\phi|_{\nu_L W}^{-1}(\tan(t))=\nu_LX_t$. Let $z\in\phi|_{\nu_L W}^{-1}(\tan(t))$ and
$z\notin \nu_LX_t$.   Using proposition 1.2, one can choose a 
subset $A$ of $W$ with $z \in \overline{A}$ and 
$\nu_L W \cap \overline{A} \cap \overline{X_t} = \varnothing$ 
(here bar denotes closure in $h_L W$).  
So for some $C, r_0>0$, 
$L(x):=\max\{ d(x, A ), d(x,X_t) \} \geq C\norm{x}$ when $x \in W$ with 
$\norm{x} \geq r_0$.  In particular, for $(y,s) \in A$ with $\norm{(y,s)} \geq r_0$, 
have 
$$ |s - \norm{y} \tan(t)| = d((y,s), (y, \norm{y}\tan(t))) \geq L(y,s) \geq C\norm{(y,s)}, $$ 
whence 
$| (s/ \norm{y}) - \tan(t) | \geq C$. But $\phi(y,s) = s / \norm{y}$.  It 
follows that $|\phi(z) - \tan(t)| \geq C$ and   
so $\phi(z) \neq \tan(t)$, a 
contradiction.     \enddemo

\proclaim{Theorem 4.3} Let $X$ be cocompact connected proper metric space.
Then $$\AN(X\times\R)=\AN X+1.$$
\endproclaim
\demo{Proof}
By Theorem 3.10, Lemma 4.2, the classical Morita theorem, and by
Theorem 3.7 we obtain
$$\AN(X\times\R)\ge\dim\nu_L(X\times\R)\ge\dim(\nu_LX\times[0,1])=
\dim\nu_LX+1\ge\AN X+1.$$
The opposite inequality is obvious.\qed
\enddemo

\head \S5 Embedding of asymptotic cones into the sublinear Higson corona \endhead

We recall the definition of the asymptotic cone $cone_{\omega}(X)$
of a metric space with the base point $x_0\in X$ with respect to a
nonprincipal ultrafilter $\omega$ on $\N$ [Gr],[Ro2]. On the
sequences of points $\{x_n\}$ with $\|x_n\|\le Cn$ for some $C$, we
define an equivalence relation
$$\{x_n\}\sim\{y_n\}\Leftrightarrow \lim_{\omega}d(x_n,y_n)/n=0.$$
We denote by $[\{x_n\}]$ the equivalence class of $\{x_n\}$. The
space $cone_{\omega}(X)$ is the set of equivalence classes
$[\{x_n\}]$ with the metric
$d_{\omega}([\{x_n\}],[\{y_n\}])=\lim_{\omega}d(x_n,y_n)/n$. We note
that the space $cone_{\omega}(X)$ does not depend on the choice of
the base point.

We note that any two constant sequences are equivalent and denote by $[x_0]$
the base point they define in $cone_{\omega}(X)$.

We call a non-principal ultrafilter $\omega$ on $\N$ {\it exponential} if it
contains the image of a function $f:\N\to \N$ satisfying the inequality
$f(n+1)\ge af(n)$ for some $a>1$ and for all but finitely many $n$.
In view of the inequality $f(n+k)\ge a^kf(n)$, the number $a$ for a given
exponential
ultrafilter $\omega$ can be taken arbitrarily large.
\proclaim{Theorem 5.1} For every exponential ultrafilter $\omega$
on $\N$ and for every proper metric space $X$, there is an injective
continuous map $\xi:cone_{\omega}(X)\setminus[x_0]\to\nu_LX$.

Moreover, the restriction of $\xi$ to the $D$-annulus
$A_D=\dot B_D([x_0])\setminus B_{1/D}([x_0])$ is an embedding for
all $D>0$.
\endproclaim
\demo{Proof}
We define
$$\xi([\{x_n\}])=\nu_LX\cap\bigcap_{S\in\omega}\overline{\{x_n\mid n\in S\}}$$
where the closure is taken in the sublinear Higson compactification $h_LX$.

First we show that for every point $[\{x_n\}]\ne[x_0]$ the above intersection is
nonempty. Indeed for every finite family of sets $S_1,\dots S_m\in\omega$
we have 
$$\nu_LX\cap\bigcap_{i=1}^m\overline{\{x_n\mid n\in S_i\}}\supset
\nu_LX\cap\overline{\{x_n\mid n\in \cap_{i=1}^kS_i\}}\ne\emptyset.$$
Then by the compactness of $\nu_LX$ we obtain the required result.

Next we show that this intersection is a point. Assume that
$x,y\in\nu_LX\cap\bigcap_{S\in\omega}\overline{\{x_n\mid n\in S\}}$, $x\ne y$.
Let $U$ and $V$ be a disjoint neighborhoods of $x$ and $y$ in $h_L X$ (in the 
sense that $U,V$ might not be open) such that 
$U \cap V = \emptyset$ and $h_LX=U \cup V$ in $h_LX$.
Let $S_x=\{n\in\N\mid x_n\in U\}$ and $S_y=\{n\in\N\mid x_n\in V\}$.
By the definition of ultrafilter, one and only one of these sets belongs
to $\omega$, say $S_x\in\omega$.
Hence $y\in\overline{\{x_n\mid n\in S_x\}}$. Therefore $y\in\bar U$ which
contradicts the fact that $V$ is a neighborhood of $y$ and 
$V$ and $U$ are disjoint.

Next we show that if $\{x_n\}\sim\{y_n\}$, then
$\xi([\{x_n\}])=\xi([\{y_n\}])$. Assume the contrary:
$$
\nu_LX\cap\bigcap_{S\in\omega}\overline{\{x_n\mid n\in S\}}\cap
\bigcap_{S\in\omega}\overline{\{y_n\mid n\in S\}}=\emptyset.$$
Then $$\nu_LX\cap\bigcap_{S\in\omega}(\overline{\{x_n\mid n\in S\}}\cap
\overline{\{y_n\mid n\in S\}})=\emptyset.$$
By compactness of $\nu_LX$ there is $\bar S\in\omega$ such that
$$
\nu_LX\cap(\overline{\{x_n\mid n\in \bar S\}}\cap
\overline{\{y_n\mid n\in \bar S\}})=\emptyset.$$
By Lemma 2.3 there is $c>0$ such that $d(x_n,\{y_k\})\ge c\|x_n\|$ for large
enough $n$. Since $[\{x_n\}]\ne [x_0]$ there is $a>0$ and $S_0\in\omega$
such that $\|x_n\|\ge an$ for $n\in S_0$. Then for $n\in S_0\cap\bar S$ we have
$d(x_n,y_n)\ge acn$. Since $\{x_n\}\sim\{y_n\}$ there is $S_1\in\omega$ such that
$d(x_n,y_n)/n<ac/2$ for $n\in S_1$. Then for $n\in S_0\cap S_1\cap \bar S$ we
obtain a contradiction: $ac/2>d(x_n,y_n)/n\ge ac$.

Next we show that $\xi$ is continuous. Let $\lim_{k\to\infty}[\{x^k_n\}]=[\{x_n\}]$.
Let $\xi([\{x_n\}])\in U$, where $U$ is an open neighborhood in $h_LX$. Then
for some $S\in\omega$ we have $\{x_n\mid n\in S\}\subset U$. By Lemma 2.3 and the fact that
$[\{x_n\}]\ne[x_0]$ there is $c>0$
such that $d(x_n,X\setminus U')\ge cn$ for some $S'\in\omega$ and some $U'$
whose closure is contained in $U$,
$\bar U'\subset U$. Let
$\lim_{\omega}d(x^k_n,x_n)/n=\delta^k$. Then $\delta^k\to 0$.
So, for large enough $k$ we have $\delta_k<c/4$.
Then there is $S_k\in\omega$
such that $d(x^k_n,x_n)/n<2\delta^k$ for $n\in S_k$. Then for $n\in S'\cap S_k$, where 
$k$ is large, 
we obtain $x^k_n\in U'$. Hence $\overline{\{x_n^k\mid n\in S'\cap S_k\}}\subset U$
for large enough $k$. Therefore $\xi([\{x^k_n\}])\in U$ for sufficiently
large $k$.

We show that $\xi$ is injective. Assume that $\xi([\{x_n\}])=\xi([\{y_n\}])$
and $[\{x_n\}]\ne[\{y_n\}]$. The latter implies that there is $\epsilon>0$ such that
$S=\{n\mid d(x_n,y_n)/n>\epsilon\}\in\omega$. We may assume
that $|\|x_n\|-D_1n|<\delta n$ and $|\|y_n\|-D_2n|<\delta n$
for $n\in S$ for some small $\delta$ where $D_1=\|[\{x_n\}]\|$
and $D_2=\|[\{y_n\}]\|$. We also may assume that $S\subset im(f)$
where $f$ is from the definition of the exponential ultrafilter
with $a\ge\max\{\frac{2D_1+D_2+3\delta}{D_2-\delta},\frac{2D_2+D_1+3\delta}{D_1-\delta}\}$. We claim that
$d(x_n,\{y_m\mid m\in S\})=d(x_n,y_n)$ and $d(y_n,\{x_m\mid m\in S\})=d(x_n,y_n)$
for $n\in S$. Indeed, for $m>n$ we have $d(x_n,y_m)\ge\|y_m\|-\|x_n\|
\ge (D_2-\delta)m-(D_1+\delta)n=(D_2-\delta)f(k+l)-(D_1+\delta)f(k)$, where 
$k$ and $l\ge 1$ are chosen such that $f(k)=n$ and $f(k+l)=m$. 
Then $d(x_n,y_m)\ge(a^l(D_2-\delta)-(D_1+\delta))f(k)\ge ((D_1+\delta)+(D_2+\delta))n\ge\|x_n\|+\|y_n\|\ge d(x_n,y_n)$. 
A similar argument works for the case $m<n$.    

By Lemma 2.3, the sets $\{x_n\mid n\in S\}$ and $\{y_n\mid n\in S\}$
diverge in the space $Z=\cup_{n\in S}\{x_n,y_n\}$ and hence in $X$.
Then 
$\nu_LX\cap\overline{\{x_n\mid n\in S\}}\cap\overline{\{y_n\mid n\in S\}}=\emptyset$, 
a contradiction.

To complete the proof it suffices to show that $\xi$
restricted to the $D$-annulus is open.

For every $[\{x_n\}]\in cone_{\omega}(X)$ with $\|[\{x_n\}]\|=c \leq D$
and for every $R$ we construct an open set $U\subset X$ such that
$\widetilde U$ contains
$\xi([\{x_n\}])$ and $\widetilde U \cap\xi(A_D)\subset\xi(B_R([\{x_n\}]))$.
Let $\lambda\ge\max\{4D/c,2cD\}$. We consider $f:\N\to \N$ with the property
$f(n+1)\ge(1+\lambda) f(n)$ and $S=Im(f)\in\omega$. Additionally we may assume
that $|\|x_n\|-cn|<\epsilon n$ for $n\in S$ for some small $\epsilon$ 
($\e < \min{ \{ c/4, 1/(4D), 1/2\} }$).  
We define $U=\cup_{n\in S} B_{\alpha\|x_n\|}(x_n)$ for
$\alpha<\min\{1/2,R/(2c)\}$.

Let $\xi([\{z_n\}])\in \widetilde U$ and $\|[\{z_n\}]\|=d$, $1/D \leq d \leq D$. 
The latter implies that
on some $S'\in\omega$,
$S'\subset S$ we have $|\|z_k\|-dk|<\epsilon k$. The former implies that
for some $S'\in\omega$ and for every $k\in S'$ there exists $n_k\in S$ such that
$d(z_k,x_{n_k})<\alpha\|x_{n_k}\|$. We may assume that in both cases we
have the same $S'$ and $S'\subset S$.

We claim that $n_k=k$.  Assume that $n_k>k$. Then the triangle inequality 
$$(d+\epsilon)k>\|z_k\|\ge\|x_{n_k}\|(1-\alpha)\ge (c-\epsilon)(1-\alpha)n_k.$$
implies that $n_k-k<\frac{4d}{c}k \leq \frac{4D}{c}k\le\lambda k$. Let $k=f(l)$. Then 
$f(l+s)=n_k$ for some $s\ge 1$ and we
obtain a contradiction: $f(l+1)-f(l)\leq f(l+s)-f(l)<\lambda f(l)$.
If we assume that $n_k<k$, then the chain of inequalities
$$(1+\alpha)(c+\epsilon)n_k>(1+\alpha)\|x_{n_k}\|\ge\|z_k\|>(d-\epsilon)k$$
implies that $k-n_k<(2c/d)n_k \leq 2cDn_k\le\lambda n_k$. If $n_k=f(l)$, then
$k=f(l+s)$ for $s\ge 1$. Then we obtain the same contradiction:
$f(l+1)-f(l)\le f(l+s)-f(l)<\lambda f(l)$.

By the construction, 
$d(z_k,x_k)<\alpha\|x_k\|\le\alpha(c+\epsilon)k \leq 2\alpha c k <Rk$.
Hence $[\{z_k\}]\in B_R([\{x_k\}])$.
\qed
\enddemo

We recall that a topological space $Y$ is called 
{\it strongly paracompact} if every open cover of $Y$ admits a star-finite
refinement.
It is known that a separable metric space is strongly paracompact
and not all metric spaces are strongly paracompact [En].

\proclaim{Corollary 5.2}
For a proper metric space $X$ and an exponential ultrafilter $\omega$, 
and for every separable $Y\subset cone_{\omega}(X)$,   
$$\dim Y\le\AN{X}.$$
\endproclaim
\demo{Proof} We present $Y=\cup_n(\bar A_n\cap Y)\cup([x_0]\cap Y)$.
Being separable metric spaces, the $\bar A_n\cap Y$ are strongly paracompact.
By 3.1.23 of [En],  $\dim Z'\le\dim Z$ for a strongly paracompact subspace $Z'\subset Z$.  
In view of Theorems 5.1 and 3.11 we obtain 
$\dim(\bar A_n\cap Y)=\dim\xi(\bar A_n\cap Y)\le\dim \nu_LX\le\AN(X)$.
The countable union theorem completes the proof.
\qed
\enddemo

QUESTION. Is it true that for a finitely generated group $G$ and an ultrafilter
$\omega$,
$$\dim cone_{\omega}G=\sup\{\dim Y\mid Y\subset cone_{\omega}G\},$$
where the supremum is taken over all separable subspaces Y?

The answer is negative if one replaces the asymptotic cone by an arbitrary metric 
space.

REMARK. Most likely Theorem 5.1 does not hold for general ultrafilters.
It seems to us that to have harmony here one should 
consider the sublinear Higson
corona $\nu_L^{\omega}X$ associated with an ultrafilter $\omega$. 
Also one can define an
Assouad-Nagata dimension $\AN_{\omega}X$ depending on $\omega$ and repeat
all basic steps of this paper to obtain the inequality 
$\dim cone_{\omega}(X)\le
\AN_{\omega}X$.

\enddocument